\documentstyle[12pt]{amsart}
\begin{document}
\begin{center}
\today\\[10pt]
{\Large\bf A Survey of Ultraproduct Constructions in General Topology} 
\\[20pt]
Paul Bankston\\
Department of Mathematics, Statistics and Computer Science\\ 
Marquette University\\
Milwaukee, WI 53201-1881\\
paulb$\mbox{@@}$mscs.mu.edu\\[20pt]

\end{center}
\begin{abstract}
We survey various attempts to transport the ultraproduct construction
from the realm of model theory to that of general topology.
\end{abstract}

\noindent
{\it A.M.S. Subject Classification\/} (1991): 03-02, 03C20, 03C52, 03C68,
54A25, 54B35, 54C10, 54D05, 54D10, 54D30, 54D80, 54E52, 54F15, 54F45,
54F50, 54F65, 54G10\\ 

\noindent
{\it Key Words and Phrases\/}:
ultraproduct, ultracoproduct, topological spaces
\\

 \noindent
{\it Dedication}: To my early mentors and colleagues:  B. Banaschewski, 
G. Bruns,
H. J. Keisler, K. Kunen, E. Nelson and M. E. Rudin.
\\[10pt]

\section{Introduction.}\label{1}
This survey is intended primarily for two readerships: general
topologists who
know a little model theory, and model theorists who know a little
general topology.  Both camps can boast a healthy constituency, and most
of this paper should be readily accessible to anyone in either.

The ultraproduct construction has a long and distinguished history.
While its beginnings go back to the 1930s with K. G\"{o}del (who was
proving his completeness theorem) and T. Skolem (who was building nonstandard
models of arithmetic), 
it was not until 1955, with the publication of the Fundamental Theorem of 
Ultraproducts, due to 
J. \L o\'{s}, that the construction was described explicitly, and its 
importance to first-order
logic became apparent. The understanding of the structure and use of 
ultraproducts developed rapidly during the next fifteen years or so,
culminating in the Ultrapower Theorem of H. J. Keisler and S. Shelah
(UT).
(The gist of the theorem is that two relational structures are elementarily
equivalent if and only if an ultrapower of one is isomorphic to an 
ultrapower of the other.  Keisler established a much stronger statement
in the early 1960s using the Generalized Continuum Hypothesis (GCH); and
toward the end of the decade, Shelah provided a GCH-free proof of a second
stronger statement that is somewhat weaker than Keisler's.)
By the late 1960s, the
theory of ultraproducts had matured into a major area of investigation
in its own right (see \cite{BS, CK, Ekl, Hod} for a vastly more detailed
account than is possible here), and was ready for export beyond the confines
of classical model theory.

Actually the exportation process had already begun by the early 1960s, when
I. Fleischer \cite{F} observed that classic ultrapowers are
directed limits of powers (and, by implication, that classic ultraproducts
are directed limits of products).  This observation, illustrating a major
strength of category theory (see \cite{MacL}), provides an abstract 
reformulation of a
concrete construction.  One may now start with a category  
{\bf C} endowed with products (which construction being itself an abstract 
reformulation of the cartesian product) 
and directed limits, and define ultraproducts within that category.   
Going further, any {\em bridging theorem\/}, i.e., one that translates a 
concrete 
notion into abstract terms involving the ultraproduct, becomes available as a 
definitional vehicle to reformulate that notion in a suitably rich category.
T. Ohkuma \cite{Ohk} (and A. Day and D. Higgs 
\cite{DH} a bit later) made good use of this idea, introducing a notion
of {\em finiteness in a category\/} by means of the important elementary 
result that says a relational
structure is finite if and only if all diagonal maps from that structure
into its ultrapowers are isomorphisms.
(In the setting of concrete categories; i.e., those endowed with a suitable
``forgetful'' functor to the category of sets and functions, this notion of 
finiteness and that of having finite underlying set can easily fail to 
coincide.  Two examples: If {\bf C} is the category {\bf CH} of 
compacta (compact
Hausdorff spaces) and continuous maps, then ``{\bf C}-finite'' means
``having at most one point.'' If {\bf C} is the category {\bf BAN} of
Banach spaces and nonexpansive homomorphisms, then ``{\bf C}-finite''
means ``being finite dimensional.'')

I became aware of Fleischer's limit approach to ultraproducts in
1974, while visiting  McMaster University, late in my career as a graduate
student.  It was there that I had the idea of using the UT as a bridging 
theorem, mimicking Ohkuma's use of the ultrapower characterization of
finiteness.  My aim was not the abstract reformulation of set-theoretic
notions, however,  but model-theoretic ones; namely elementary equivalence and
elementary embedding (as well as their various derivative notions). 
I can attribute much of my own development as a mathematician
to enlightening talks I had with the universal algebra group at McMaster
at that time (namely B. Banaschewski, G. Bruns and E. Nelson), and the
papers \cite{Ban6,Ban7,BF} extend and develop the ideas introduced in
\cite{Ohk, DH}.  Moreover, my coinage of the term ``ultracoproduct,'' along  
with my own investigations of how ultraproducts behave in the {\em opposite\/} 
of the concrete category {\bf CH} (to be discussed in \S 5) 
can also be traced to Fleischer's approach. 

What Fleischer started in 1963 might be regarded as the beginning of
the idea of a ``model-theoretic study of  
a class (or category) {\bf C}.''  This should be immediately contrasted 
with what might be called
``{\bf C}-based model theory.''  While the two subject areas may overlap a 
great deal, there is a difference in emphasis.  In the former, one perhaps
fixes an autonomous notion of ultraproduct in {\bf C} (hence a mechanism for 
generating conjectures that stem from known classical results),
then tries to establish (functorial) links between {\bf C}
and particular classes of models of first-order theories (hence a mechanism
for settling some of those conjectures).  In the latter,
one enriches objects of {\bf C} with
extra ``functions'' and ``relations,'' possibly nonclassical in nature but
recognizable nonetheless, views these enriched objects as models of logical 
languages, and proceeds to develop new model theories, using more established
model theories for guidance.
Our study of compacta in \cite{Ban12} and elsewhere exemplifies the 
former emphasis, while the Banach model theory initiated by C. W. Henson
(see \cite{Hen1,Hen2, HI}), as well as the approaches to topological model 
theory
found in \cite{FZ, Gar, McK, Sgro}, exemplify the latter.

In this paper, our primary focus is on how classical ultraproducts can be
exported to purely topological contexts, with or without category-theoretic
considerations as motivation.
(So the Banach ultraproduct
\cite{DK}, for example, the Fleischer ultraproduct in {\bf BAN}, is not
directly a subject of our survey.)  We begin in the next section 
with a quick introduction to ultraproducts in model theory; then on, in 
\S 3, to consider the topological ultraproduct, the most straightforward
and na\"{\i}ve attempt at exporting the ultraproduct to the topological context.
The motivation in \S 3 is purely model-theoretic, with no overt use of
category-theoretic concepts.
This is also true in \S 4, where we look at a variation of this construction 
in the special case of ultrapowers.  It is not until \S 5, where 
ultracoproducts are introduced, that the Fleischer approach to defining
ultraproducts plays a significant role.  Although the ultracoproduct may
be described in purely concrete (i.e., set-theoretic) terms, and is of
independent interest as a topological construction, the important 
point is that category-theoretic language allows one to see this construction 
as a natural gateway out of the classical model-theoretic context.

The ultraproduct construction in model theory is a quotient of the direct
product, where an ultrafilter on the index set dictates how to specify the 
identification.  When we carry out the analogous process in general
topology, at least from the viewpoint of \S 3, the ``product'' in
question is not the usual Tychonov product, but the less sophisticated
(and much worse-behaved) {\it box\/} product.  (While one could use the
Tychonov product instead of the box product, the result would be 
an indiscrete (i.e., trivial) topological space, unless the ultrafilter
were countably complete.)   

The identification process just mentioned does not require the maximality
of the designated ultrafilter in order to be well defined, 
and may still be carried out using any filter 
on the index set.  The resulting construction, called the 
{\it reduced product\/}, serves as a generalization of both the direct 
(box) product and the ultraproduct constructions.  In \S 6 we survey some
of the recent work on furthering this generalization to include the
Tychonov product and some of its relatives.  Finally, in \S 7, we list
some of the more resistant and intriguing open problems in the topological
study of ultraproducts.\\  

\section{Preliminaries from Model Theory.}\label{2}
First we recall some familiar notions from model theory, establishing our 
basic notation and terminology in the process.

Given a set $I$, the power set of $I$ is
denoted ${\wp}(I)$, and is viewed as a bounded lattice under unions and
intersections. (The alphabet of bounded lattices consists of two binary
operation symbols, $\sqcup$ (join) and $\sqcap$ (meet), plus two constant
symbols, $\top$ (top) and $\bot$ (bottom).) 
A {\bf filter on} $I$ is a filter {\em in} the lattice
${\wp}(I)$; i.e., a collection $\cal F$ of subsets of $I$ satisfying:
$(i)$ $I \in {\cal F}$, $(ii)$ any superset of an element of $\cal F$ is
also an element of $\cal F$, and $(iii)$ the intersection of any two elements 
of $\cal F$ is also an element of $\cal F$.  A filter $\cal F$ is called 
{\bf proper} if $\emptyset \notin {\cal F}$; an {\bf ultrafilter on} $I$
is a proper filter on $I$ that is not contained in any other (distinct)
proper filter on $I$; i.e., a maximal proper filter in the lattice 
${\wp}(I)$.  In power set lattices, the maximal proper filters are
precisely the {\bf prime} ones; that is, any proper filter $\cal F$ on $I$
is an ultrafilter, if for each $J,K \subseteq I$, if $J \cup K \in {\cal F}$,
then either $J \in {\cal F}$ or $K \in {\cal F}$.
If ${\cal S}$ is any family of subsets of $I$, 
$\cal S$ is said to satisfy the {\bf Finite Intersection Property
(FIP)} if no finite intersection of elements of $\cal S$ is empty.  Our
underlying set theory of choice is Zermelo-Fraenkel Set Theory 
{\em with} Choice
(ZFC); consequently, any family of subsets of $I$ that satisfies the FIP
must be contained in an ultrafilter on $I$.  (More generally, if a subset
of a bounded distributive lattice satisfies the Finite {\em Meet} Property,
then that subset is contained in a maximal proper filter in the lattice.)

We start with an alphabet $L$ of
finitary relation and function symbols (with the equality symbol $\approx$
tacitly assumed to be included).  
An $L$-{\bf structure}
consists of an underlying set $A$ and an interpretation of each symbol
of $L$, in the usual way.  Like many authors (and unlike many others), we
use the same font to indicate both a relational structure and its underlying
set; being careful to make the distinction clear whenever there is a threat
of ambiguity.

If $\langle A_i:i\in I\rangle$ is an indexed family of $L$-structures, 
and $\cal F$ is a filter on $I$, the ordinary direct product of the family is 
denoted $\prod_{i\in I}A_i$, with the $i$th co\"{o}rdinate of an element
$a$ being denoted $a(i)$.  (Each symbol of $L$ is interpreted in the  
standard way.)  The binary relation $\sim_{\cal F}$ on the
product, given by $a \sim_{\cal F} b$ just in case $\{i\in I: a(i)=b(i)\}
\in {\cal F}$, is easily seen to be an equivalence relation; and we
define $a/{\cal F} := \{b : a \sim_{\cal F} b\}$.  We denote by
$\prod_{\cal F}A_i$ the corresponding {\bf reduced product}; i.e., the set of 
$\sim_{\cal F}$-equivalence classes, with the standard interpretation of
each symbol of $L$.  When $A_i = A$ for each $i \in I$, we have the
{\bf reduced power}, denoted $A^I/{\cal F}$.  The {\bf canonical diagonal map} 
$d: A \to A^I/{\cal F}$, given by $a \mapsto (\mbox{constantly}\;a)/{\cal F}$,
is clearly an embedding of $L$-structures.  
  
From here on, unless we specify otherwise, we concentrate
on reduced products (powers) in which the filter is an ultrafilter.  The
corresponding constructions are called {\bf ultraproducts (ultrapowers)},
and the Fundamental Theorem of Ultraproducts is the following. (We follow
the standard notation regarding satisfaction of substitution instances of
first-order formulas.  That is, if $\varphi (x_0,\dots,x_{n-1})$ is a 
first-order $L$-formula with free variables from the
set $\{x_0,\dots,x_{n-1}\}$, and if $A$ is an $L$-structure with
$n$-tuple $\langle a_0,\dots,a_{n-1}\rangle \in A^n$, then $A \models
\varphi[a_0,\dots,a_{n-1}]$ means that the sentence got from $\varphi$
by substituting each free occurrence of $x_i$ with a new constant
symbol denoting $a_i$, $i<n$, is true in $A$.  (See also
\cite{BS,CK,Hod}.))

\subsection{Theorem.}\label{2.1} 
(\L o\'{s}' Fundamental Theorem of Ultraproducts \cite{CK}) 
Let $\langle A_i: i\in I\rangle$ 
be a family of $L$-structures, with $\cal D$ an ultrafilter on $I$ and
$\varphi (x_0,\dots,x_{n-1})$ a first-order $L$-formula.  Given an
$n$-tuple $\langle a_0/{\cal D},\dots,a_{n-1}/{\cal D}\rangle$ from
the ultraproduct, then $\prod_{\cal D}A_i \models
\varphi[a_0/{\cal D},\dots,a_{n-1}/{\cal D}]$ if and only if
$\{i\in I: A_i\models \varphi[a_0(i),\dots,a_{n-1}(i)]\}\in {\cal D}$.\\

By a {\bf level zero formula}, we mean a Boolean combination of atomic
formulas.  If $k$ is any natural number, define a {\bf level $k+1$
formula} to be a level $k$ formula $\varphi$ preceded by a string $Q$ of 
quantifiers of like parity (i.e., either all universal or all existential) 
such that, if $\varphi$ begins with a quantifier, then the parity of
that quantifier is not the parity of the quantifiers of $Q$.
Formulas with a well-defined level are said to be in {\bf prenex form}, and
elementary first-order logic provides
an effective procedure for converting any $L$-formula to a logically 
equivalent formula (with the same free variables) in prenex form.
A function $f:A \to B$ between $L$-structures is a
{\bf level $\geq k$ embedding} if for each $L$-formula 
$\varphi (x_0,\dots ,x_{n-1})$ of level $k$, 
and $n$-tuple $\langle a_0,\dots ,a_{n-1}\rangle
\in A^n$, it is the case that $A \models \varphi[a_0,\dots ,a_{n-1}]$ if
and only if $B\models \varphi[f(a_0),\dots ,f(a_{n-1})]$.
It is easy to see that the level $\geq 0$ embeddings are precisely the 
algebraic embeddings;
the level $\geq 1$ embeddings are also called {\bf existential embeddings}.
(Existential embeddings have been of considerable interest to algebraists
and model theorists alike.)  If a function $f$ is of level $\geq k$ for all
$k < \omega$, we call it a {\bf level $\geq \omega$ embedding}.  Now an
{\bf elementary embedding} is one that preserves the truth of all first-order
formulas, even those without an obvious level; so elementary embeddings are
clearly of level $\geq \omega$.  The effective procedure mentioned above,
then, assures us of the converse.  We are taking pains to make this point
because, as we shall see, the notion of level $\geq k$ embedding can be
given a precise abstract meaning, devoid of reference to first-order
formulas.  Moreover, it can be extended into the transfinite, giving rise
to an abstract notion of level $\geq \alpha$ {\em morphism\/.}   There is
no {\em a priori\/} reason that this hierarchy should terminate at level
$\omega$. (See, e.g., \S 5.)
    
\subsection{Corollary.}\label{2.2}
(Diagonal Theorem) The canonical diagonal embedding from a relational
structure into an ultrapower of that structure is an elementary embedding.\\

A first-order formula containing no free variables is called a {\bf sentence},
and two $L$-structures $A$ and $B$ are called {\bf elementarily equivalent} 
(denoted $A \equiv B$) if they
satisfy the same $L$-sentences.
Clearly if there is an elementary
embedding from one $L$-structure into another, then the two structures 
are elementary equivalent; in particular, because of \ref{2.2}, if
some ultrapower of $A$ is isomorphic to some ultrapower of $B$, then
$A \equiv B$.  By the same token, if $f:A \to B$ is a map between 
$L$-structures, then $f$ is an elementary embedding as long as there are
ultrafilters $\cal D$ and $\cal E$ (on sets $I$ and $J$ respectively) and an 
isomorphism $h:A^I/{\cal D} \to B^J/{\cal E}$ such that the compositions
$e\circ f$ and $h\circ d$ with the corresponding diagonal embeddings are equal. 
The converses of these two statements are also true. (Indeed, the converse
of the second follows from the converse of the first via the method of
expanding the alphabet $L$ by adding constants denoting all the elements of
$A$.)  This fact, called the (Keisler-Shelah) Ultrapower Theorem (UT), is a
milestone in model theory, with a very interesting history
(see, e.g., \cite{CK}).  Its importance, in part, is that it allows many basic
notions of first-order model theory to be formulated in abstract
terms, i.e., in terms of mapping diagrams; it is what we called a bridging  
theorem in the Introduction.  The obvious central notions are elementary
equivalence and elementary embedding, but there are also derivative
notions (e.g., prime model) readily definable in terms of these.  Other
derivative notions are less obvious.  The following result is stated
and used extensively in \cite{HSim}, and is an application of Keisler's 
Model Extension Theorem (see \cite{Sho}).  

\subsection{Theorem.}\label{2.3}
A function $f:A \to B$ between $L$-structures is a level $\geq k+1$ embedding 
if and only if 
there is an elementary embedding $e:A \to C$ and a level $\geq k$ embedding
$g:B \to C$ such that $e = g\circ f$.\\

Theorem \ref{2.3}, in conjunction with the UT, is another bridging
theorem; as the elementary embedding $e$ may be taken to be an
ultrapower diagonal embedding.  Thus the notion of level $\geq k$
embedding has an abstract reformulation.  Indeed, because of the inductive
flavor of \ref{2.3}, that notion may be formally carried into the
realm of transfinite levels.  What is more, the notion of level $\geq 1$
(existential) embedding is now available in abstract form.  That means
we can export model completeness to the category-theoretic setting. 

We begin to see how these ideas may be exploited when we survey the
topological ultracoproduct in \S 5.  (We use the infix {\em co\/} because
we are dealing with the {\em opposite} of the concrete category {\bf CH}.) 
Using the UT as a bridge, abstract model-theoretic notions are imported,
only in dual form, and made concrete once again.  In order for this to
be a productive enterprise, however, it is necessary to
use more of the theorem than simply the ``gist'' form stated above.  
We therefore
end this section with statements of both
Keisler's GCH version and Shelah's subsequent GCH-free version.
(We employ standard notation as regards cardinals and ordinals; see,
e.g., \cite{CK}.  In particular, if $\kappa$ and $\lambda$ are cardinals,
then $\kappa^+$ is the cardinal successor of $\kappa$; and $\kappa^{\lambda}$
is the cardinal exponential, the cardinality of the set of all functions from
$\lambda$ into $\kappa$.  If $S$ is any set, its cardinality is denoted
$|S|$.)

\subsection{Theorem.}\label{2.4}
(Keisler's Ultrapower Theorem, \cite{CK,She}) Let $\lambda$ be an infinite 
cardinal where
the GCH holds (i.e., $2^{\lambda} = \lambda^+$), and let $I$ be a set
whose cardinality is $\lambda$.  Then there is an ultrafilter
$\cal D$ on $I$ such that if $L$ is an alphabet with at most $\lambda$ 
symbols, and if $A$ and $B$ are elementarily equivalent
$L$-structures of cardinality at most $\lambda^+$, then $A^I/{\cal D}
\cong B^I/{\cal D}$.\\ 

\subsection{Theorem.}\label{2.5}
(Shelah's Ultrapower Theorem. \cite{She}) Let $\lambda$ be an infinite 
cardinal, with $\mu := \mbox{min}\{\alpha: \lambda^{\alpha}>\lambda\}$,
and let $I$ be a set whose cardinality is $\lambda$.  Then there is an
ultrafilter $\cal D$ on $I$ such that if $L$ is an alphabet, and if
$A$ and $B$ are elementarily equivalent $L$-structures of cardinality
less than $\mu$, then $A^I/{\cal D} \cong B^I/{\cal D}$.

\section{Topological Ultraproducts.}\label{3}
Following established usage, a {\bf topological space} consists of an
underlying set $X$ and a family $\cal T$ of subsets of $X$, called a 
{\bf topology}; members of $\cal T$ being called {\bf open sets}.  
All a family of subsets has to do to be called a topology
is to be closed under arbitrary unions and finite intersections.  As with
the case of relational structures, we use the same symbol to indicate both a
topological space and its underlying set (using disambiguating notation, 
such as $\langle X,{\cal T}\rangle$, only when necessary).  If $\cal B$ is
an open base for a topology $\cal T$ on $X$ (so arbitrary unions of members
of $\cal B$ form a topology), then we write ${\cal T} = \tau ({\cal B})$,
the topology {\bf generated by} $\cal B$. 

Let $\langle \langle X_i,{\cal T}_i\rangle : i\in I\rangle$ be an indexed 
family of topological
spaces, with $\cal D$ an ultrafilter on $I$.  Then the ultraproduct
$\prod_{\cal D}{\cal T}_i$ may easily be identified with a family of subsets
of the ultraproduct $\prod_{\cal D}X_i$, and this family qualifies as an
open base for a topology $\tau ({\prod_{\cal D}{\cal T}_i})$
on $\prod_{\cal D}X_i$, which we call the 
{\bf ultraproduct topology}.  The resulting {\bf topological ultraproduct}
is denoted (when we can get away with it) $\prod_{\cal D}X_i$; and the
canonical basic open sets $\prod_{\cal D}U_i \in \prod_{\cal D}{\cal T}_i$
are called {\bf open ultraboxes}.  Clearly the quotient map $x \mapsto
x/{\cal D}$ from $\prod_{i\in I}X_i$ to $\prod_{\cal D}X_i$ is a continuous
open map from the box product to the ultraproduct.   
Also there is a certain amount of flexibility built into the definition of
topological ultraproduct; in that one may obtain an open base for the
ultraproduct topology by taking open ultraboxes $\prod_{\cal D}U_i$, where,
for each $i \in I$, the sets $U_i$ range over an open base for the
topology ${\cal T}_i$.  This flexibility extends to {\it closed\/} bases
as well.  Recall that a family $\cal C$ is a {\bf closed base} for
$\cal T$ if $\cal T$-closed sets (i.e., complements in $X$ of members of
$\cal T$) are intersections of subfamilies taken from $\cal C$.  One may
obtain a closed base for the ultraproduct topology by taking {\it closed\/}
ultraboxes $\prod_{\cal D}C_i$, where, for each $i \in I$, the sets
$C_i$ range over a closed base for the topology ${\cal T}_i$.
(The reader interested in nonstandard topology may want to compare the
topological ultrapower topology with A. Robinson's {\it $Q$-topology\/}
\cite{Rob}.)

The connection between topological ultraproducts and usual ultraproducts
should be rather apparent, but we will find it convenient to spell things
out.  By the {\bf basoid alphabet} we mean the alphabet $L_{BAS} :=
\{P,B,\varepsilon\}$, where the first two symbols are unary
relation symbols standing for ``points'' and ``basic open sets,'' 
respectively, and the third, a binary relation symbol, stands for 
``membership.''  If $X$ is any set and ${\cal S} \subseteq {\wp}(X)$,
then $\langle X,{\cal S}\rangle$ may be naturally viewed as the
$L_{BAS}$-structure 
$\langle X \cup {\cal S},X,{\cal S},\in\rangle$,
where set-theoretic membership is restricted to $X\times {\cal S}$. 
An $L_{BAS}$-structure is called a {\bf basoid} if it is (isomorphic to)
such a structure, where $\cal S$ is an open base for a topology on $X$.
The basoid is called {\bf topological} if $\cal S$ is itself a topology.
Every basoid has a uniquely associated topological basoid; the second
is said to be {\bf generated} from the first.
It is a routine exercise to show that
there is a first-order $L_{BAS}$-sentence whose models are
precisely the basoids.  Thus ultraproducts of basoids are basoids by
Theorem \ref{2.1}, and we obtain 
$\prod_{\cal D}\langle X_i,{\cal T_i}\rangle$ as the topological basoid
generated from the usual ultraproduct of the basoids 
$\langle X_i,{\cal T_i}\rangle$. 

The alphabet $L_{BAS}$ is a natural springboard for topological model
theory:  Allow extra relation and function symbols to range over points,
and build various languages from there.  This is a one-sorted approach,
which is quite sensible, but which turns out to be somewhat cumbersome in 
practice for the purposes of exposition.  Other approaches in the 
literature start with a first-order alphabet $L$, and expand the first-order
language over $L$ in various ways.  For example, there is 
the extra-quantifiers approach, exemplified by
J. Sgro's $L_Q$ \cite{Sgro} (patterned after Keisler's $L_Q$ \cite{Kei});
also the two-sorted approach, exemplified by the ``invariant'' languages
$L_t$ of T.A. McKee \cite{McK} and S. Garavaglia \cite{Gar}.  (The two
worked independently, with McKee confining himself to the case $L = 
\{\approx\}$.  See also \cite{FZ}.)  There is an extensive  model theory
for $L_t$ which we cannot possibly survey adequately.  (The 
interested reader is urged to consult
the Flum-Ziegler monograph \cite{FZ}.)  However, since this model theory
includes a nice ultrapower theorem, 
we take a few lines to describe these languages and state the theorem.

One starts with an ordinary first-order
alphabet $L$, adds new variables to stand for sets, and then adds the
intersorted binary relation symbol $\varepsilon$ for membership. 
Atomic formulas consist of the first-order atomic formulas from $L$,
plus the intersorted formulas of the form $t \varepsilon U$, where $t$
is a first-order term (from $L$) and $U$ is a set variable.  The language
$L_2$ consists of the closure of the atomic formulas under the
logical connectives $\neg$ (not), $\vee$ (or) and $\wedge$ (and), 
and the quantifiers $\exists$ (there exists) and 
$\forall$ (for all), applied to variables of either sort.  
A formula $\varphi$ of $L_2$ is {\bf positive} (resp., {\bf negative}) in
the set variable $U$ if each free occurrence of $U$ in $\varphi$ lies
within the scope of an even (resp., odd) number of negation symbols.
We then define $L_t$ to be the smallest subset $K$ of $L_2$ satisfying:
$(i)$ the atomic formulas are in $K$; $(ii)$ $K$ is closed under the logical
connectives, as well as quantification over point variables; and $(iii)$ if
$t$ is a first-order term and $\varphi \in K$ is positive (resp., negative)
in $U$, then
$(\forall U(\neg(t \varepsilon U) \vee \varphi)) \in K$
(resp., $(\exists U((t \varepsilon U) \wedge \varphi)) \in K$).

By a {\bf basoid $L$-structure}, we mean a pair $\langle A, {\cal B}\rangle$,
where $A$ is an $L$-structure and $\cal B$ is an open base for some topology
on $A$.  It should then be clear what it means for a basoid structure to
be a {\it model\/} of a sentence $\varphi$ of $L_2$, as well as what it means
for two basoid $L$-structures to be {\it isomorphic}.  If 
$\langle A_1,{\cal B}_1\rangle$ and
$\langle A_2,{\cal B}_2\rangle$ are two basoid $L$-structures, then these
structures are {\bf homeomorphic} just in case 
$\langle A_1,\tau({\cal B}_1)\rangle$ and
$\langle A_2,\tau({\cal B}_2)\rangle$ are isomorphic.  We may now state
the topological version of the ultrapower theorem, due to
Garavaglia, as follows. 

\subsection{Theorem.}\label{3.1}
(Garavaglia's Ultrapower Theorem \cite{FZ, Gar}) Let $A$ and $B$ be two
basoid $L$-structures.  Then $A$ and $B$ satisfy the same $L_t$-sentences
if and only if some ultrapower of $A$ is homeomorphic to some ultrapower
of $B$.\\

In \cite{Ban2}, two spaces $X$ and $Y$ are said to be {\bf power equivalent}
if some ultrapower of $X$ is homeomorphic to some
ultrapower of $Y$.  It is not hard to show directly (Theorem A2.3 in 
\cite{Ban2}) that power equivalence is really an equivalence relation, and
it is of some interest to see just how strong an equivalence relation it
is.  Recall that a space is said to be {\bf self-dense} if it has no isolated 
points.
We use the well-known $T_n$-numbering of the separation axioms
(\`{a} la \cite{Wil}); but note that, for the purposes of this paper, 
we assume the $T_1$ axiom (i.e., singletons are closed) whenever we talk
about separation axioms involving arbitrary closed sets.  Thus
{\it regularity\/} (resp., {\it normality\/}), the property
of being able to separate a point and a non-containing closed set 
(resp., two disjoint closed sets) with
disjoint open sets, presupposes the $T_1$ axiom, and is synonymous with
the $T_3$ (resp., $T_4$) axiom.  Similarly, we assume $T_1$ when we define
{\it complete regularity\/} (or, the {\it Tychonov property\/}, sometimes
referred to as the $T_{3.5}$ axiom) as the 
property of being
able to separate a point and a non-containing closed set with a continuous
real-valued function.  
The following tells us that
power equivalence is not very discriminating.

\subsection{Theorem.}\label{3.2}
(Theorem A2.6 in \cite{Ban2}) Any two self-dense $T_3$-spaces are power
equivalent.\\

\subsection{Remark.}\label{3.3}
The proof of \ref{3.2} uses a combination of model theory and topology.
In particular, it makes use of the L\"{o}wenheim-Skolem theorem and
a result of W. Sierpi\'{n}ski \cite{Sier}, to the effect that any two countable,
second countable, self-dense $T_3$-spaces are homeomorphic.  One could
claim that \ref{3.2} is a corollary of \ref{3.1}, but that would
be a stretch.  One would still need to employ the theorems of 
L\"{o}wenheim-Skolem and Sierpi\'{n}ski to show that any two 
topological basoids that are self-dense and $T_3$ must satisfy the same
$L_t$-sentences (where $L = \{\approx \}$).\\

With any apparatus that produces new objects from old, an important
issue concerns the idea of preservation.  In the context of the topological
ultraproduct construction, a preservation problem takes the following
general form.

\subsection{Problem.}\label{3.4}
(General Preservation) Given topological properties $P$ and $Q$, and  
a property $R$ of ultrafilters, 
decide the following:  For any $I$-indexed family $\langle X_i:i \in I\rangle$
of topological spaces and any ultrafilter $\cal D$ on $I$, if 
$\{i\in I: X_i\;\mbox{has property}\;P\}\in {\cal D}$ 
(i.e., ``$\cal D$-almost every $X_i$ has property $P$'')
and $\cal D$ has property $R$, then $\prod_{\cal D}X_i$ has
property $Q$.\\

\subsection{Remark.}\label{3.5}
The general problem, as stated in \ref{3.4}, is not quite as general as
it could be.  The property $P$ could actually be a family $\cal P$ of
properties, and the clause 
``$\cal D$-almost every $X_i$ has property $P$'' could read 
``$\cal D$-almost every $X_i$ has property $P$ for all  
$P \in {\cal P}$.'' The vast majority of instances of this
problem do not require the added generality, however. (One obvious exception: 
Consider, for $n < \omega$, the property $P_n$ that 
says that there are at least $n$ points, and set ${\cal P}:= \{P_n:
n<\omega\}$.  If $R$ is the property of being countably incomplete
and $Q$ is the property of being infinite, then this instance of 
the more general version of \ref{3.4}
has an affirmative answer.) \\ 

The question of the preservation 
of the separation axioms $T_0$--$T_4$ under ultraproducts turns out to
be a very rich topic.
In \cite{Ban2} we define a topological property $P$ to be {\bf closed}
if \ref{3.4} has an affirmative answer for $Q=P$ and $R$ nonrestrictive.
$P$ is {\bf open} if its negation is closed.  It is a straightforward
exercise in definition manipulation to show \cite{Ban2} that the 
axioms $T_0$ through $T_3$ are both closed and open; a little less 
straightforward to show is the fact that $T_{3.5}$ is closed.  
It should come as no surprise to general topologists that neither $T_{3.5}$
nor $T_4$ is open, and that $T_4$ is not closed.  The proofs of these
negative facts are fairly involved, but they are
valuable more for how they themselves involve an assortment of new ideas and
positive (general) results.  For this reason  
we take a few paragraphs to expatiate on some of their key points.

Consider first why $T_{3.5}$ is a closed property.  Recall the
well-known characterization of O. Frink \cite{Wal} that a $T_1$-space $X$ is
completely regular if and only if 
it has a {\bf normal disjunctive lattice base}; that is, if there is
a bounded sublattice $\cal C$ of the bounded lattice of closed subsets
of $X$ satisfying: $(i)$ $\cal C$ is a closed base for the topology on $X$
(i.e., $\cal C$ is meet-dense in the closed set lattice); $(ii)$ (normality)
for each disjoint pair $C,D \in {\cal C}$ there exist $C',D' \in {\cal C}$
with $C \cap C' = D \cap D' = \emptyset$ and $C' \cup D' = X$; and 
$(iii)$ (disjunctivity) for each two distinct elements of $\cal C$, there
is a nonempty element of $\cal C$ that is contained in one of the first
two elements and is disjoint from the other. (A good source on basic
distributive lattice theory is \cite{Balb}.)  
If $\langle X_i: i \in I\rangle$
is a family of spaces such that $\cal D$-almost every $X_i$ is completely
regular, then for $\cal D$-almost every $i \in I$, there is a normal
disjunctive lattice base ${\cal C}_i$ for $X_i$.  It follows quickly that
$\prod_{\cal D}{\cal C}_i$ is a normal disjunctive lattice base for
$\prod_{\cal D}X_i$.

Of the twelve preservation results above concerning $T_0$--$T_4$, only the 
first nine are apparently
positive.  Nevertheless, it so happens that the last three are corollaries  
of positive results.  Indeed,  
one can show that both $T_{3.5}$ and $T_4$ are not open properties in
one go, with the help of \ref{3.2}.  A space $X$ is {\bf linearly orderable}
(a {\bf LOTS}) if $X$ has a linear ordering whose open intervals constitute an
open base for $X$.  $X$ is {\bf linearly uniformizable} (a {\bf LUTS}) if    
the topology on $X$ is induced by a uniformity that has a linearly ordered
base under inclusion.  (See, e.g., \cite{Wil}.  For example, if $\rho$
is a metric inducing the topology on $X$, then $\{\{\langle x,y\rangle:
\rho(x,y) < \epsilon\}:\epsilon > 0\}$ is a linearly ordered uniform base
that witnesses the fact that $X$ is a LUTS.)  Suppose $\cal D$-almost
every $X_i$ is a LOTS with inducing linear order $\leq_i$ (resp., a LUTS with
inducing linearly ordered uniform base ${\cal U}_i$).  Then 
$\prod_{\cal D}\leq_i$ (resp., $\prod_{\cal D}{\cal U}_i$) is a linear
ordering (resp., a linearly ordered uniform base) that induces the
ultraproduct topology on $\prod_{\cal D}X_i$.  
Now every LOTS is hereditarily
normal; indeed every LUTS is hereditarily paracompact Hausdorff.  So let $X$ 
be any regular space.  Then $X\times \Bbb R$, the topological product of $X$
with the real line, is self-dense and regular.  By \ref{3.2}, there
is an ultrapower $(X\times \Bbb R)^I/{\cal D}$ that is homeomorphic
to an ultrapower of $\Bbb R$, and is hence both a LOTS and a LUTS.
It is easy to show that ultrapowers commute with finite products.  
Thus the ultrapower $X^I/{\cal D}$ embeds in an ultrapower of the 
reals, and is hence hereditarily normal (indeed, hereditarily paracompact
Hausdorff).  
The following theorem, whose proof we have just outlined, immediately
implies the failure of $T_{3.5}$ and $T_4$ to be open properties.

\subsection{Theorem.}\label{3.6}
(Corollary A2.7 in \cite{Ban2}) Every regular space has a hereditarily 
paracompact Hausdorff ultrapower.\\

We now turn to the problem of showing that normality is not a closed property.
First some notation:
If $\kappa$ and $\lambda$ are cardinals, we write $\kappa^{\lambda}$
to indicate the $\lambda$-fold topological (Tychonov) power of the ordinal
space $\kappa$ (as well as the cardinal exponentiation).  The following
positive result clearly implies that normality fails to be closed.  

\subsection{Theorem.}\label{3.7}
(Corollary of Theorem 8.2 in \cite{Ban2}) Let $X$ be any space
that contains an embedded copy of $2^{\omega_2}$, and let $\cal D$ be 
any nonprincipal ultrafilter on a countable set $I$.  Then $X^I/{\cal D}$
is not normal.\\

The proof of \ref{3.7}, being far more interesting than the statement, 
deserves a bit of discussion.

Of course, if $X$ fails to be regular, so does any ultrapower.  Thus
it suffices to confine our attention to regular $X$ (or even Hausdorff;
it does not matter).  In that case, any embedded copy $Y$ of $2^{\omega_2}$
is closed in $X$; hence $Y^I/{\cal D}$ is closed in $X^I/{\cal D}$.
It is therefore enough to show that $Y^I/{\cal D}$ is nonnormal.

This brings us to the important class of $P$-spaces.  
Following the terminology of \cite{GJ, Wal}, we
call a space $X$ a {\bf $P$-space} if every countable intersection of
open sets is an open set.  More generally, following the Comfort-Negrepontis 
text \cite{CN}, let $\kappa$
be an infinite cardinal.  A point $x$ in a space $X$ is called a 
{\bf $P_{\kappa}$-point} if for every family $\cal U$ of fewer than $\kappa$ 
open neighborhoods of $x$, there is an open neighborhood of $x$ that
is contained in each member of $\cal U$.  $X$ is a {\bf $P_{\kappa}$-space}
if each point of $X$ is a $P_{\kappa}$-point.  In $P_{\kappa}$-spaces,
intersections of fewer than $\kappa$ open sets are open; the $P$-spaces
are just the $P_{\omega_1}$-spaces.  (In \cite{Ban1,Ban2}, the 
$P_{\kappa}$-spaces
are called $\kappa$-{\it open}.  While it is convenient to have a concise 
adjectival form of ``being a $P_{\kappa}$-space,'' there was already one
in the literature, $\kappa$-{\it additive\/}, due to R. Sikorski \cite{Sik}, 
which we adopt here.)

It is very hard for a topological ultraproduct not to be a $P$-space.  To
be specific, define an ultrafilter $\cal D$ on $I$ to be {\bf $\kappa$-regular}
if there is a family ${\cal E} \subseteq {\cal D}$, of cardinality $\kappa$,
such that each member of $I$ is contained in only finitely many members
of $\cal E$.  It is well known \cite{CN} that $|I|^+$-regular ultrafilters 
cannot exist, that $|I|$-regular ultrafilters exist in abundance, 
that $\omega$-regularity is the same as countable incompleteness, and
that nonprincipal ultrafilters on countable sets are countably incomplete.
The following not only says that $\kappa$-regularity in ultrafilters produces
$\kappa^+$-additivity in topological ultraproducts (deciding affirmatively
an instance of Problem \ref{3.4}); it actually 
{\it characterizes\/} this property of ultrafilters. 

\subsection{Theorem.}\label{3.8}   
(Additivity Lemma, Theorem 4.1 in \cite{Ban2}) An ultrafilter is 
$\kappa$-regular if and only if all topological ultraproducts via that
ultrafilter are $\kappa^+$-additive ($P_{\kappa^+}$-spaces).\\
 
\subsection{Remark.}\label{3.9}
There is a model-theoretic analogue to \ref{3.8}:  Just replace ``additive''
with ``universal.'' (See Theorem 4.3.12 and Exercise 4.3.32 in \cite{CK}.)\\

Given any space $X$ and cardinal $\kappa$, we denote by $(X)_{\kappa}$
the space whose underlying set is $X$, and whose topology is the smallest
$\kappa$-additive topology containing the original topology of $X$.  If
$\kappa$ is a regular cardinal (so $\kappa$ is not the supremum of fewer
than $\kappa$ smaller cardinals; for example $\kappa$ could be a successor
cardinal), then one may obtain an open base for
$(X)_{\kappa}$ by taking intersections of fewer than $\kappa$ open subsets
of $X$.  (See, e.g., \cite{CN} for an extensive treatment of this kind of
topological operation.) 

\subsection{Terminological Remark.}\label{3.10}
The adjective {\it regular\/}, as used in technical mathematics, is probably
the most overloaded word in the English language.  Already in this paper
it has three senses; modifying the nouns {\it space\/},
{\it ultrafilter\/} and {\it cardinal\/} in completely unrelated ways.   
In other areas of mathematics as well, the word is used with abandon.
In algebra, functions, 
rings, semigroups, permutations and representations can all be regular;
in homotopy theory, fibrations can be regular;
and in analysis, Banach spaces, measures and points can be regular too.  
({\it Regular\/} modifies {\it ring\/} in the same way that it modifies
{\it semigroup\/}, but otherwise there are no apparent similarities in the 
senses to which it is used.)  The list, I am sure, goes on.\\

Returning to the proof outline of \ref{3.7}, recall the diagonal map $d$
from a set $X$ into an ultrapower $X^I/{\cal D}$ of that set.  If the
ultrapower is a topological one, $d$ is not necessarily continuous; 
consider, for example the case where $X$ is the real line and $\cal D$
is a countably incomplete ultrafilter.  The image $d[X]$ of $X$ under $d$ then
carries the discrete topology.  The following uses \ref{3.8}.

\subsection{Theorem.}\label{3.11} 
(Theorem 7.2 in \cite{Ban2}) Let $\cal D$ be a regular ultrafilter on a
set of cardinality $\kappa$, with $X$ a topological space.  Then the diagonal
map, as a map from $(X)_{\kappa^+}$ to $X^I/{\cal D}$, is a topological
embedding.\\

Suppose $Y$ is a compactum (i.e., a compact Hausdorff space), and that
$\cal D$ is an ultrafilter on $I$.  Then for each $a/{\cal D} \in
X^I/{\cal D}$, there is a unique point $x \in X$ such that for each open
set $U$ containing $x$, the open ultra{\it cube\/} $U^I/{\cal D}$ contains
$a/{\cal D}$.  Let $\lim_{\cal D}(a/{\cal D})$ denote this unique point.     
Then the function $\lim_{\cal D}$ is continuous (Theorem 7.1 in \cite{Ban2}),
and is related to the standard part map in nonstandard analysis \cite{Rob}.  
But more is true, thanks to \ref{3.11}.

\subsection{Theorem.}\label{3.12}
(A consequence of Corollary 7.3 of \cite{Ban2})
Let $\cal D$ be a regular ultrafilter on a set of cardinality $\kappa$,
with $Y$ a compactum.  Then the limit map $\lim_{\cal D}$, as a map
from $Y^I/{\cal D}$ to $(Y)_{\kappa^+}$, is a continuous left inverse
for the diagonal map $d$.  As a result, the diagonal $d[Y]$, a homeomorphic
copy of $(Y)_{\kappa^+}$, is a closed subset of $Y^I/{\cal D}$.\\

We are just about done with \ref{3.7}.  In a preliminary version of 
\cite{Kun},  K. Kunen shows that
$(2^{{\frak c}^+})_{{\omega}_1}$ is nonnormal, where ${\frak c} := 
2^{{\aleph}_0}$ is the power of the continuum; and in \cite{vanD},
E. K. van Douwen uses an earlier result of C. Borges \cite{Bor} to replace 
$\frak c$
with $\omega_1$.  So let $Y$ now be the compactum $2^{{\omega}_2}$, with
$\cal D$ any nonprincipal ultrafilter on a countable set $I$.  In order
to show $Y^I/{\cal D}$ is nonnormal, it suffices to show some closed
subset is nonnormal.  This is true, though, since $(Y)_{\omega_1}$ is
nonnormal and, by \ref{3.12}, sits as a closed subset of $Y^I/{\cal D}$. 
This completes our discussion of \ref{3.7}.
 
What Borges' result cited above actually says is that the space
$({\kappa}^{\kappa^+})_{\kappa}$ is nonnormal whenever $\kappa$ is a
regular cardinal.  It is quite easy to show from this that, for any
infinite cardinal $\kappa$, $(2^{{\kappa}^{++}})_{\kappa^+}$ is not
normal either.  This, together with the additivity lemma \ref{3.8} and
some arguments to show how easy it is for paracompactness to be present
in $P$-spaces, gives rise to a characterization of the GCH in terms of
topological ultraproducts.  

Recall that the {\bf weight} of a space $X$ is the greater of $\aleph_0$
and the least cardinality of an open base for the topology on $X$.
For each infinite cardinal $\kappa$, let UP$_{\kappa}$ be the following
assertion.\\ 

\noindent
UP$_{\kappa}$: If $I$ is a set of cardinality $\kappa$, $\cal D$ is a
regular ultrafilter on $I$, and $\langle X_i:i \in I\rangle$ is an
$I$-indexed family of spaces, $\cal D$-almost each of which is regular
and of weight at most $2^{\kappa}$, then $\prod_{\cal D}X_i$ is
paracompact Hausdorff.\\

The main result of \cite{Ban1} (see also W. Comfort's survey article
\cite{Com}) is the following.

\subsection{Theorem.}\label{3.13}
(Theorem 1.1 in \cite{Ban1}) UP$_{\kappa}$ holds if and only if
the GCH holds at level $\kappa$ (i.e., $2^{\kappa} = {\kappa}^+$).\\

\subsection{Remarks.}\label{3.14}
$(i)$ The proof of \ref{3.13} allows several alternatives to UP$_{\kappa}$.
In particular, {\it regular\/} (as the word applies to spaces) may be
replaced by {\it normal\/}; even by {\it compact Hausdorff\/}.  Also
{\it paracompact Hausdorff\/} may be replaced by {\it normal\/}. 

$(ii)$ Topological ultraproducts are continuous open images of box
products, and there are many inevitable comparisons to be made between
the two constructions.  In particular, let BP$_{\kappa}$ be the statement
that the box product of a $\kappa$-indexed family of compact Hausdorff
spaces, each of weight at most $2^{\kappa}$, is paracompact Hausdorff.  
In \cite{Kun} it is proved that the CH (i.e., the GCH at level $\omega$)
implies BP$_{\omega}$.  Since $(2^{\omega_2})_{\omega_1}$ is nonnormal,
the compactum $2^{\omega_2}$ stands as a counterexample to 
BP$_{\omega}$ if the CH fails, and as an {\em absolute} counterexample to
BP$_{\kappa}$ for $\kappa > \omega$.\\

We now turn to the exhibition of Baire-like properties in topological
ultraproducts.
If $\kappa$ is an infinite cardinal, define a space $X$ to be
{\bf $\kappa$-Baire} (or, a {\bf $B_{\kappa}$-space}) if intersections
of fewer than $\kappa$ dense open subsets of $X$ are dense.  Of course,
every space is a $B_{\omega}$-space, and various forms of the Baire
category theorem say that completely metrizable spaces and compact
Hausdorff spaces are $\omega_1$-Baire.  Finally, one topological form of
Martin's Axiom (MA, see, e.g., \cite{BMc}) says that if $X$ is compact
Hausdorff and satisfies the {\it countable chain condition\/} (i.e.,
there is no uncountable family of pairwise disjoint nonempty open subsets
of $X$), then $X$ is $\frak c$-Baire.  

What we are working toward is an analogue of \ref{3.8}, with $P$
replaced with $B$.  What has been achieved in this connection
is interesting, if imperfect, and begs for improvement.

For any set $S$ and cardinal $\lambda$, let ${\wp}_{\lambda}(S)$ be the
set of all subsets of $S$ of cardinality less than $\lambda$. 
If $\cal D$ is an ultrafilter on a set $I$, a map $F:{\wp}_{\omega}(S)
\to {\cal D}$ is {\bf monotone} (resp., {\bf multiplicative}) if $F(s)\supseteq
F(t)$ whenever $s \subseteq t$ (resp., $F(s \cup t) = F(s) \cap F(t)$).  
The ultrafilter $\cal D$ is called {\bf $\lambda$-good} if: $(i)$ $\cal D$
is countably incomplete, and $(ii)$ for every $\mu < \lambda$ and every
monotone $F: {\wp}_{\omega}(\mu) \to {\cal D}$, there exists a
multiplicative 
$G: {\wp}_{\omega}(\mu) \to {\cal D}$ such that $G(s) \subseteq F(s)$
for all $s \in {\wp}_{\omega}(\mu)$.  (This notion is due to Keisler.) 

Every countably incomplete ultrafilter is $\omega_1$-good, and every
$\lambda$-good ultrafilter is $\mu$-regular for all $\mu < \lambda$.
Consequently, if $|I| = \kappa$, the maximal degree of goodness an 
ultrafilter on $I$ could hope to have is $\kappa^+$.  The existence
of good ultrafilters (i.e., $\kappa^+$-good ultrafilters on sets of
cardinality $\kappa$) was first proved by Keisler under the hypothesis
$2^{\kappa} = \kappa^+$, and later by Kunen without this hypothesis.
(See \cite{CN}.  There it is shown that there are as many good ultrafilters
on a set as there are ultrafilters.)  Good ultrafilters produce saturated  
models (see Theorem 6.1.8 in \cite{CK}), and the production of saturated
models necessitates goodness (see Exercise 6.1.17 in \cite{CK}).  Finally, 
and most importantly, good ultrafilters play a crucial role in the proofs
of both ultrapower theorems \ref{2.4} and \ref{2.5}.  Our analogue of
\ref{3.8} is the following affirmative answer to the general preservation 
problem (\ref{3.4}).   

\subsection{Theorem.}\label{3.15}
(Theorem 2.2 in \cite{Ban3}) If an ultrafilter is $\kappa$-good, then all
topological ultraproducts via that ultrafilter are $\kappa$-Baire
($B_{\kappa}$-spaces) (as well as being $\lambda^+$-additive for all
$\lambda < \kappa$).\\
 
\subsection{Remark.}\label{3.16}
Theorem \ref{3.8} is actually key to the proof of \ref{3.15}.  We do not
know whether producing topological ultraproducts that are $\kappa$-Baire as 
well as $\lambda^+$-additive for all $\lambda < \kappa$ is sufficient to show
an ultrafilter to be $\kappa$-good.\\

Topological ultraproduct methods have proven useful in the study of
the $\eta_{\alpha}$-sets of F. Hausdorff \cite{Hau}.  
Recall that, for any infinite cardinal $\alpha$, 
a linear ordering $\langle A,<\rangle$
is an {\bf $\eta_{\alpha}$-set} if whenever $B,C \subseteq A$ each have
cardinality less than $\alpha$, and every element of $B$ lies to the
left of every element of $C$, then there is some element of $A$ lying
to the right of every element of $B$ and to the left of every element of $C$.
The $\eta_{\omega}$-sets are just the dense linear orderings without endpoints,
and Hausdorff \cite{Hau} invented the famous ``back and forth'' method to
show that any two $\eta_{\alpha}$-sets of
cardinality $\alpha$ are order isomorphic.  He was also able to establish
the existence of $\eta_{\alpha^+}$-sets of cardinality $2^{\alpha}$ (and
L. Gillman showed how to exhibit two distinct such orderings whenever 
$\alpha^+ <
2^{\alpha}$).  Gillman and B. J\'{o}nsson proved that  
$\eta_{\alpha}$-sets of cardinality $\alpha$ exist precisely under the 
condition that $\alpha = \sup\{\alpha^{\lambda}:\lambda < \alpha\}$.
(The interested reader should consult \cite{GJ,CN}.)
Denote by ${\Bbb Q}_{\alpha}$ the (unique, when it exists) 
$\eta_{\alpha}$-set of cardinality $\alpha$.  (${\Bbb Q}_{\omega}$ is,
of course, the rational line $\Bbb Q$.)  In \cite{Ban9}, we use 
topological ultraproduct methods to establish properties of 
${\Bbb Q}_{\alpha}$, viewed as a LOTS.  In particular, 
${\Bbb Q}_{\alpha}$ is both $\alpha$-additive and $\alpha$-Baire, and the
following is true.

\subsection{Theorem.}\label{3.17}  
(Theorem 3.14 of \cite{Ban9}) If $X$ is a nonempty space that embeds in
${\Bbb Q}_{\alpha}$, then  
${\Bbb Q}_{\alpha}$ can be partitioned into homeomorphic copies of $X$,
each of which is closed and nowhere dense in 
${\Bbb Q}_{\alpha}$.\\ 

We end this section with one more preservation result about topological
ultraproducts.  Its main interest is that its proof apparently needs to
involve two cases, depending upon whether the ultrafilter is countably
complete or countably incomplete.  Also it involves a topological property 
that illustrates a general machinery for producing new properties from old. 

By \ref{3.8}, every topological ultraproduct via a countably incomplete
ultrafilter is a $P$-space.  Now if a $P$-space is also $T_1$, then
it has the peculiar property of being {\bf pseudofinite} (or,
a {\bf $cf$-space}, see \cite{HR}); i.e., one 
having no infinite compact subsets.  Another way of saying this is that
the only compact subsets of $X$
are the ones that have to be, based on cardinality considerations alone.

There is a general phenomenon afoot here.  Namely,
if $P$ is any topological property, let $\mbox{spec}(P)$
be the set of cardinals $\kappa$ such that every space of cardinality
$\kappa$ has property $P$; and denote by $\mbox{anti-}P$ the class of
spaces $X$ such that if $Y$ is a subspace of $X$ and $Y$ has property
$P$, then $|Y| \in \mbox{spec}(P)$.  For example, if $P$ is the property
{\it compact\/} (resp., {\it connected\/}, {\it self-dense\/}), then 
$\mbox{anti-}P$ is the property {\it pseudofinite\/}
(resp., {\it totally disconnected\/}, {\it scattered\/}). 
The modifier {\it anti-\/} was introduced in \cite{Ban5}, and it has
been studied in its own right by a number of workers.  (See, e.g.,
\cite{RV,MMc1,MMc2}.)  Concerning topological ultraproducts, what we showed 
in \cite{Ban5} is the following affirmative answer to \ref{3.4}.

\subsection{Theorem.}\label{3.18}
(Corollary 3.6 of \cite{Ban5}) Topological ultraproducts of pseudofinite
Hausdorff spaces are pseudofinite Hausdorff.\\

\subsection{Remark.}\label{3.19}
Of course, topological ultraproducts of Hausdorff spaces, via countably
incomplete ultrafilters, are pseudofinite Hausdorff (by \ref{3.8} plus
basic facts).  One must argue quite differently when the ultrafilters
are countably complete.  In this case cardinal measurability is involved, 
and pseudofiniteness on the part of
the factor spaces is essential; moreover the argument does not work if
the Hausdorff condition is eliminated (or even weakened to $T_1$).  
One needs to know that if a set has a certain cardinality, then the
cardinality of its closure cannot be too much greater.  The $T_2$ axiom
assures us of this, but the $T_1$ axiom does not.  (Consider any set
with the cofinite topology.)  So, for example, 
we do not know whether topological ultraproducts of pseudofinite 
$T_1$-spaces are pseudofinite in general.\\

\section{Coarse Topological Ultrapowers.}\label{4}
There is a natural variation on the definition of the ultraproduct  
topology in cases where all the factor spaces are the same. 
In this section, we consider ultrapowers only, and restrict the ultrapower
topology to the one generated by just the open ultracubes.  This is what
we call the {\bf coarse topological ultrapower}.  That is, if 
$\langle X,{\cal T}\rangle$ is a topological space and $\cal D$ is an 
ultrafilter on a set $I$, then the family of open ultracubes $\{U^I/{\cal D}:
U \in {\cal T}\}$ forms an open base for the coarse ultrapower topology.
Note that, with regard to this topology, the natural diagonal map
$d: X \to X^I/{\cal D}$ is a topological embedding.   
We denote the coarse topological ultrapower by $[X^I/{\cal D}]$.
(For those interested in nonstandard topology, there is a connection
between coarse topological ultrapowers and Robinson's {\it $S$-topology\/}
\cite{Rob}.)

Quite straightforwardly, one may obtain a closed base for the coarse
ultrapower topology by taking all closed ultracubes.  However, it is 
generally {\it not} true that an open (resp., closed) base for the
coarse ultrapower topology may be obtained by taking ultracubes from
an open (resp., closed) base for the original space.  (Indeed, let
$X$ be infinite discrete, with ${\cal B}$ the open base of singleton
subsets of $X$.)  

Our main interest in this section is the question of
when coarse topological ultrapowers satisfy any of the usual separation
axioms.  If the ultrafilter is countably complete, then the diagonal map
is a homeomorphism unless the base space has cardinality exceeding the
first measurable.  While this may be an interesting avenue of research,
there are no results at this time that we know of; and we therefore 
confine attention to countably incomplete ultrafilters.  For each   
$r \in \{0,1,2,3,3.5,4\}$, define an
ultrafilter $\cal D$ to be a {\bf $T_r$-ultrafilter} if it is countably
incomplete, and for some infinite space $X$, the coarse ultrapower
$[X^I/{\cal D}]$ is a $T_r$-space.  The reader should have no difficulty
in constructing coarse topological ultrapowers that are not $T_0$-spaces,
so the question of the mere existence of $T_0$-ultrafilters will doubtless
come to mind.  The good news is that $T_0$-ultrafilters are closely
related to ones whose combinatorial properties are fairly well understood; so
their existence follows from MA.  We currently do not know whether
$T_0$-ultrafilters exist absolutely.    

First, we may reduce the existence question to the case of ultrafilters
on a countable set; $\omega$, say.  The reason is that if $\cal D$ is
a $T_0$-ultrafilter on an infinite set $I$ and $X^I/{\cal D}$ is $T_0$,
then we may partition $I$ into countably many subsets, none of which
is in $\cal D$, and build a function $f$ from $I$ onto $\omega$ such that
the images of the members of the partition of $I$ partition $\omega$ into
infinite sets.  Then ${\cal E} := \{S \subseteq \omega: f^{-1}[S] \in
{\cal D}\}$ is clearly a countably incomplete ultrafilter.  Moreover
$f$ induces an embedding of $X^{\omega}/{\cal E}$ into $X^I/{\cal D}$;
hence $\cal E$ is a $T_0$-ultrafilter. 

In \cite{Sco}, B. Scott defines an ultrafilter $\cal D$ on $\omega$ to
be {\bf separative} if whenever $f,g:\omega \to \omega$ are two functions
that are {\bf $\cal D$-distinct} (i.e., $\{n<\omega: f(n) \neq g(n)\}
\in {\cal D}$), then their Stone-\v{C}ech lifts $f^{\beta}$ and $g^{\beta}$
disagree at the point ${\cal D} \in \beta(\omega)$ (i.e., there is some
$J \in {\cal D}$ such that $f[J] \cap g[J] = \emptyset$).  Scott's main
results in \cite{Sco} include the facts that selective
ultrafilters are separative, and the properties of selectivity and being
a $P$-point ultrafilter (i.e., a $P$-point in $\beta (\omega)\setminus \omega$)
are not implicationally related.  From MA, one may infer the
existence of selective ultrafilters; hence the consistency of separative  
ultrafilters is assured.  By the famous Shelah $P$-point independence
theorem \cite{Wim}, $P$-point ultrafilters cannot be shown to exist in
ZFC.  We do not know whether the same can be said for separative ultrafilters,
but strongly suspect so.  The following is an amalgam of several results in
\cite{Ban8}.

\subsection{Theorem.}\label{4.1}
An ultrafilter on $\omega$ is $T_r$, for $r \in \{0,1,2,3,3.5\}$, if and
only if it is separative.\\

\subsection{Remark.}\label{4.2}
That $\cal C$ is separative if it is $T_0$ is straightforward (Proposition
2.1 in \cite{Ban8}).  Assuming $\cal D$ is separative, it is shown in
\cite{Ban8} that a coarse $\cal D$-ultrapower of $X$ is: $(i)$ $T_1$ if
$X$ is a {\it weak $P$-space\/} (i.e., no point is in the closure
of any countable subset of the complement of the point) ; $(ii)$ $T_2$ if
$X$ is $T_2$ and a $P$-space; $(iii)$ $T_{3.5}$ if $X$ is $T_4$
and a weak $P$-space; and $(iv)$ {\it strongly zero dimensional\/}
(i.e., disjoint zero sets are separable via disjoint closed open sets) if
$X$ is $T_4$ and a $P$-space.  We do not know whether coarse topological
ultrapowers (of infinite spaces, via countably incomplete ultrafilters)
can {\it ever\/} be normal.\\  

\section{Topological Ultracoproducts.}\label{5}
Most algebraists at all familiar with the classical reduced product construction
know how to define it in terms of direct limits of products (\`{a} la
Fleischer \cite{F}).  Indeed, in his introductory article in the
``Handbook of Mathematical Logic,'' Paul Eklof \cite{Ekl} goes this
route, but then says:

\begin{quote}
``Although the shortest approach to the definition of reduced products is via
the notion of direct limit, this approach is perhaps misleading since it is 
the concrete construction of the direct limit rather than its universal
mapping properties which will be of importance in the sequel.''
\end{quote}

Eklof quite sensibly proceeds immediately to the concrete construction
(i.e., in terms of elements), because classical model theory has no use for
the abstract approach.  However, there is more to ultraproducts than
just first-order logic.  Consider, for example, the problem of giving an
explicit concrete description of the Stone space of an ultraproduct of
Boolean lattices, in terms of the Stone spaces of those lattices.  
(Note:  We speak of Boolean {\em lattices\/}, rather than Boolean 
{\em algebras\/}, because we do not include complementation as a distinguished
unary operation.  Of course, in the context of bounded distributive lattices,
complements are unique when they exist.  Thus a bounded sublattice of a
Boolean algebra is a subalgebra just in case the sublattice is Boolean 
itself.)  
Because of the duality theorem of M. H. Stone (see
\cite{Joh}), this ``ultracoproduct'' must be an inverse limit of coproducts.
To be more definite, suppose $\langle X_i:i \in I\rangle$ is an $I$-indexed
family of Boolean (i.e., totally disconnected compact Hausdorff) spaces, 
with $\cal D$ an ultrafilter on $I$.  Letting $B(X)$ denote the Boolean 
lattice of ``clopen'' (i.e., closed open) subsets of $X$, the operator $B(\;)$
is contravariantly functorial, with ``inverse'' given by the maximal
spectrum functor $S(\;)$.  Given any Boolean lattice 
$A$, the points of $S(A)$ are
the maximal proper filters in $A$.  If $a \in A$ and $a^{\sharp} :=
\{M \in S(A): a \in M\}$, then the set $A^{\sharp}:= \{a^{\sharp}: a \in A\}$
forms a (closed) lattice base for a totally disconnected
compact Hausdorff topology on $S(X)$.

So Stone Duality tells us that $S(\prod_{\cal D}B(X_i))$ is an inverse
limit of coproducts; hence a subspace of $\beta(\bigsqcup_{i\in I}X_i)$,
the Stone-\v{C}ech compactification of the disjoint union of the spaces
$X_i$.  Here is one way (out of many) to describe this space in purely
topological terms. Let $Y$ be   
$\bigsqcup_{i\in I}X_i$, and let $q:Y \to I$ take an element to its index.
Then there is the natural Stone-\v{C}ech lift $q^{\beta}:\beta(Y) \to 
\beta(I)$ ($I$ having the discrete topology), and it is not hard to show 
that $S(\prod_{\cal D}B(X_i))$ is naturally homeomorphic to 
$(q^{\beta})^{-1}[\cal D]$, the inverse image of ${\cal D} \in \beta(I)$ under
$q^{\beta}$.  Let us denote this space $\sum_{\cal D}X_i$.  It is rightfully
called an ``ultracoproduct'' because it is category-theoretically dual to
the usual ultraproduct in a very explicit way.  What makes this whole
exercise interesting is that our explicit description of 
$\sum_{\cal D}X_i$ requires nothing special about the spaces $X_i$ 
beyond the Tychonov separation axiom.  Indeed, the construction just
described, what we call the {\bf topological ultracoproduct}, is the 
Fleischer-style ultraproduct for the {\em opposite} of the category {\bf CH}
of compacta and continuous maps.  And while the topological ultracoproduct
makes sense in the Tychonov context ($\sum_{\cal D}X_i$ is actually a
compactification of the topological ultraproduct $\prod_{\cal D}X_i$), 
one does not get anything new in
the more general setting.  That is, 
$\sum_{\cal D}X_i$ is naturally homeomorphic to $\sum_{\cal D}\beta (X_i)$
(see \cite{Ban12}).  For this reason we confine our attention to 
ultracoproducts of compacta.  

If each $X_i$ is the same compactum $X$, then we have the
{\bf topological ultracopower}
$XI\backslash {\cal D}$, a
subspace of $\beta(X\times I)$. In this case there is the Stone-\v{C}ech
lifting $p^{\beta}$ of the natural first-co\"{o}rdinate map 
$p:X\times I\to X$.  Its restriction to the ultracopower is a continuous
surjection, called the {\bf codiagonal map}, and is officially 
denoted $p_{X,{\cal D}}$ (with the occasional notation-shortening alias
possible). 
This map is dual to the natural diagonal map from a relational structure
to an ultrapower of that structure, and is 
not unlike the standard part map from nonstandard analysis.  (It is closely
related to, indeed an extension of, $\lim_{\cal D}$, introduced after
\ref{3.11}.)  Recalling 
from the Introduction that being abstractly finite in {\bf CH} means
having at most one point (because abstract ultraproducts via countably
incomplete ultrafilters must have trivial topologies),  it is natural
to ask what being abstractly finite means in the opposite of {\bf CH}.
The answer is simple, satisfying and easy to prove: it means 
``having a finite number of points.''

Stone Duality is a contravariant equivalence between the categories 
{\bf BS} of Boolean spaces and continuous maps and {\bf BL} of 
Boolean lattices and homomorphisms.  From our perspective, {\bf BL} is
an interesting participant in the duality because it has abstract products, 
all cartesian, and its class of objects is one that
is first-order definable.  This tells us its Fleischer-style ultraproduct
construction is the usual one.   
For the purposes of this paper, let us
call a concrete category {\bf C} {\bf Stone-like} if there is a contravariant
equivalence between {\bf C} and some concrete category {\bf A}, with usual
(cartesian) products; where the objects of
{\bf A} are the models of a first-order theory, and the morphisms of 
{\bf A} are the functions that preserve atomic formulas.  Then clearly any
Stone-like category has an ultracoproduct construction, in the Fleischer
sense of forming inverse limits of coproducts.  Thus {\bf BS} is
Stone-like, as is the category {\bf CAG} of compact Hausdorff abelian
groups and continuous group homomorphisms.  The reason {\bf CAG} is 
Stone-like is that there is a celebrated duality theorem, due to L. Pontryagin
(see, e.g., \cite{Wil}), that matches this category with the category {\bf AG}
of abelian groups and homomorphisms.  But while the ultraproduct constructions
in {\bf BL} and {\bf AG} are exactly the same, the ultracoproduct constructions
in {\bf BS} and {\bf CAG} are quite different \cite{Ban14}. 
 
Any time a concrete category {\bf C} has an abstract ultra(co)product
construction, there are two clear lines of investigation that present
themselves.  First one may study the construction {\em per se\/} in 
set-theoretic 
terms, via the underlying set functor; second one may view the construction 
as a vehicle for establishing abstract formulations of various model-theoretic
notions (thanks to the Ultrapower Theorem).  The second line is more
``global'' in flavor; it is part of a study of the category {\bf C} as
a whole.  For example, one may wish to know whether {\bf C} is Stone-like.
(As explained in \cite{Ban6}, the full subcategory {\bf TDCAG} of totally
disconnected compact Hausdorff abelian groups, a category with an abstract
ultracoproduct construction, is not Stone-like because it has ``cofinite''
objects with infinite endomorphism sets.)
As one might expect, it is
a combination of these two lines that gives the best results. 
  
Now we have seen that there is an abstract ultraproduct construction, as
well as an abstract ultracoproduct construction, in the category {\bf CH}.
As we saw earlier, the first construction is uninteresting because it
almost always has the trivial topology.  The story is quite different
for the second, however.  For one thing, it extends the corresponding
construction in the full subcategory {\bf BS}, so there is an immediate
connection with model-theoretic ultraproducts.  (In fact there is generally
a natural isomorphism between $B(\sum_{\cal D}X_i)$ and 
$\prod_{\cal D}B(X_i)$.  This implies, of course, that ultracoproducts
of connected compacta (=continua) are connected \cite{Ban6, Ban12}.)  
For another thing, there is the fact that a compactum $X$ is finite if and only
if all codiagonal maps $p_{X,{\cal D}}$ are homeomorphisms (``cofinite''
= finite).  

In light of the above, a natural conjecture to make is that {\bf CH} is
Stone-like; and after over twenty years, everything known so far about
the topological ultracoproduct points to an affirmative answer (in
contrast to the situation with {\bf TDCAG}).
I first posed the question in the McMaster algebra seminar in 1974, and
expressed then my belief that the conjecture is false.  At the time I had
little more to go on than the empirical observation that there were already
quite a few duality theorems involving {\bf CH}, e.g.,
those of Banaschewski, Morita, Gel'fand-Kolmogorov and Gel'fand-Na\u{\i}mark,
and none of them were of the right kind. 
Almost ten years (and several partial answers, see \cite{Ban6}) later, 
there came confirmation 
of my belief from two independent quarters.

\subsection{Theorem.}\label{5.1}
(B. Banaschewski \cite{Bana} and J. Rosick\'{y} \cite{Ros}) {\bf CH} is
not a Stone-like category.\\

Of course, what Banaschewski and Rosick\'{y} independently prove are two
somewhat different-sounding statements that each imply \ref{5.1}.  
The importance of their result is that it underscores the point that
dualized model-theoretic analogues of classical results, automatically
theorems in Stone-like categories, are merely conjectures in {\bf CH}.
(Shining example:  R. L. Vaught's Elementary Chains Theorem.)
  
Because of the failure of {\bf CH} to be Stone-like (perhaps this 
``failure'' is a virtue in disguise), one is forced to look elsewhere
for model-theoretic aids for a reasonable study of topological ultracoproducts.
Fortunately there is a finitely axiomatizable universal-existential Horn class 
of bounded distributive lattices, 
the so-called {\bf normal disjunctive} lattices (also called
Wallman lattices ), comprising precisely the (isomorphic
copies of) {\bf lattice bases}, those lattices that serve as bases for the 
closed sets of compacta.  (To be more specific:  The normal disjunctive
lattices are precisely those bounded lattices $A$ such that there exists
a compactum $X$ and a meet-dense sublattice ${\cal A}$ of the closed set
lattice $F(X)$ of $X$ such that $A$ is isomorphic to $\cal A$.)  
We go from bounded distributive lattices to spaces, as in the case of Stone
duality, via the {\bf maximal spectrum} $S(\;)$, pioneered by H. Wallman
\cite{Walm}. 
$S(A)$ is the space of 
maximal proper filters of $A$; a typical basic closed set in $S(A)$ is the
set $a^{\sharp}$ of elements of $S(A)$ containing a given element $a \in A$.  
$S(A)$ is generally compact with this topology.
Normality, the condition that if
$a$ and $b$ are disjoint ($a \sqcap b = \bot$), then there are $a'$,
$b'$ such that $a \sqcap a' = b \sqcap b' = \bot$ and $a' \sqcup b' = \top$,
ensures that the maximal spectrum topology is Hausdorff.  Disjunctivity,
which says that for any two distinct lattice elements there is a nonbottom 
element that is below one
and disjoint from the other,
ensures that the map $a \mapsto a^{\sharp}$ takes
$A$ isomorphically onto the canonical closed set base for $S(A)$.  $S(\;)$
is contravariantly functorial: If $f:A \to B$ is a homomorphism of normal
disjunctive lattices and $M \in S(B)$, then $f^S(M)$ is the unique maximal
filter extending the prime filter $f^{-1}[M]$.  (For normal 
lattices, each prime filter is contained in a unique maximal one.)

It is a relatively easy
task to show, then, that $S(\;)$ converts ultraproducts to ultracoproducts.
Furthermore, if $f:A\to B$ is a {\bf separative} embedding; i.e.,
an embedding such that
if $b \sqcap c = \bot$ in $B$, then there exists $a \in A$ such that
$f(a) \geq b$ and $f(a) \sqcap c = \bot$, then $f^S$ is a homeomorphism. 
Because of this, there is much
flexibility in how we may obtain $\sum_{\cal D}X_i$:  Simply choose a
lattice base ${\cal A}_i$ for each $X_i$ and apply $S(\;)$ to the
ultraproduct $\prod_{\cal D}{\cal A}_i$.  So, taking each ${\cal A}_i$
to be $F(X_i)$, we infer very quickly that $\sum_{\cal D}X_i$ contains
the topological ultraproduct $\prod_{\cal D}X_i$ (\`{a} la \S 3) as
a densely-embedded subspace.  
(Also we get an easy concrete
description of the codiagonal map $p:XI\backslash {\cal D} \to X$:  If
$\cal A$ is a lattice base for $X$ and $y \in XI\backslash {\cal D}
= S({\cal A}^I/{\cal D})$, then $p(y)$ is that unique $x \in X$ such
that if $A \in {\cal A}$ contains $x$ in its interior, then 
$A^I/{\cal D} \in y$.  So $p$ does indeed extend $\lim_{\cal D}$.  
Note that we may view the $\cal D$-equivalence class
of the ``constantly $A$'' function in ${\cal A}^I$ as an ultrapower
itself.)    

So we officially define two compacta $X$ and $Y$ to be {\bf co-elementarily
equivalent} if there are ultracopowers $p: XI\backslash {\cal D} \to X$,
$q: YJ\backslash {\cal E} \to Y$, and a homeomorphism $h: XI\backslash {\cal D}
\to YJ \backslash {\cal E}$. (Recall the definition of {\em power equivalence\/}
in \S 3.)  A function $f:X \to Y$ is a {\bf co-elementary
map} if there are $p$, $q$, and $h$ as above such that the compositions
$f\circ p$ and $q\circ h$ are equal.  These definitions come directly from
the UT.  Furthermore, because of Theorem \ref{2.3}, we may define the
{\bf level} of a map $f:X \to Y$ as follows:  $f$ is a map {\bf of level
$\geq 0$} if $f$ is a continuous surjection.  If $\alpha$ is any ordinal,
$f$ is a map {\bf of level $\geq \alpha +1$} if there are maps
$g:Z \to Y$ and $h:Z \to X$ such that $g$ is co-elementary, $h$ is of 
level $\geq \alpha$, and $f\circ h = g$.  If $\alpha$ is a positive limit
ordinal, $f$ is a map {\bf of level $\geq \alpha$} if $f$ is a map of level
$\geq \beta$ for all $\beta < \alpha$.  (Because of the definition of 
{\em co-elementary map\/}, $g: Z \to Y$ may be taken to be an ultracopower
co-diagonal map.)  A map of level $\geq 1$ is also called {\bf
co-existential}.   

The reader may be wondering whether we are justified in the terminology
``co-elementary {\em equivalence\/},'' as there is nothing in the definition
above that ensures the transitivity of this relation.  The answer is that
we are so justified; but we need the maximal spectrum functor $S(\;)$, plus 
the full power of the UT (i.e., Theorem \ref{2.5}) to show it 
(Theorem 3.2.1 in \cite{Ban12}).  By the same token, one also shows that
compositions of co-elementary maps are co-elementary (Theorem 3.3.2 in
\cite{Ban12}), and that compositions of maps of level $\geq \alpha$ are
of level $\geq \alpha$ (Proposition 2.5 in \cite{Ban20}).   

Because of how it translates ultraproducts of lattices to ultracoproducts of
compacta, the maximal spectrum functor also translates elementary 
equivalence between lattices to co-elementary equivalence between compacta.
Furthermore, if $f:A \to B$ is an elementary (resp. level $\geq \alpha$)
embedding, then $f^S:S(B) \to S(A)$ is a co-elementary (resp. level
$\geq \alpha$) map.  Nevertheless, 
the spectrum functor falls far short of being a duality, except when restricted
to the Boolean lattices.  For this reason, one must take care not to jump
to too many optimistic conclusions; such as assuming, e.g., that if  
$f:X \to Y$ is a co-existential map, then there must
be lattice bases $\cal A$ for $X$ and $\cal B$ for $Y$ and an existential
embedding $g: {\cal B} \to {\cal A}$ such that $f = g^S$.  
(Of course,
for level $\geq 0$, this is obvious:  Pick ${\cal A} := F(X)$, ${\cal B} := 
F(Y)$,
and $g := f^F$.  However, $f^F$ is not an existential embedding, unless it is 
already an isomorphism (a slight adjustment of the proof of Proposition 2.8 
in \cite{Ban18}).) This 
``representation problem'' has yet to be solved.

The infrastructure for carrying out a dualized model-theoretic study of
compacta is now in place.  Because of Stone Duality, dualized model theory
for Boolean spaces is perfectly reflected in the ordinary model theory
of Boolean lattices, but \ref{5.1} tells us there is no hope for a similar
phenomenon in the wider context.  For example, one may use the Tarski
Invariants Theorem \cite{CK}, plus Stone Duality, to show that there are
exactly $\aleph_0$ co-elementary equivalence classes in {\bf BS}; however,
one must work directly to get the number of co-elementary equivalence
classes in {\bf CH}.

\subsection{Theorem.}\label{5.2}
(Diversity, Theorem 3.2.5 in \cite{Ban12}) 
There are exactly $\frak c$ co-elementary
equivalence classes in {\bf CH}. Moreover (Theorem 1.5 in \cite{Ban14}),
for each $0 < \alpha \leq \omega$, there is a family of $\frak c$ metrizable 
compacta, each of dimension $\alpha$, no two of which are co-elementarily
equivalent.  
Finally (Theorem 2.11 in \cite{Ban15}),
there is a family of $\frak c$ locally connected metrizable (i.e., {\em
Peano\/}) continua, no two of which are co-elementarily equivalent.\\ 

Another example concerns various statements of the L\"{o}wenheim-Skolem
Theorem.  The weakest form, for Boolean lattices, says that every Boolean
lattice is elementarily equivalent to a countable one (``countable'' =
``countably infinite or finite'').  Now Stone Duality equates the cardinality
of a Boolean lattice with the weight of its maximal spectrum space
(in symbols, $|A| = w(S(A))$); hence
we infer immediately that every Boolean space is co-elementarily equivalent
to a metrizable one (since, for compacta, metrizability = weight $\aleph_0$).
The same is true for compacta in general, by use of the 
L\"{o}wenheim-Skolem Theorem for normal disjunctive lattices.  This was 
first proved by R. Gurevi\v{c} \cite{Gur}, in response to a question raised
in \cite{Ban12}. 

\subsection{Theorem.}\label{5.3}
(L\"{o}wenheim-Skolem Theorem, Proposition 16 in \cite{Gur}) For every
compactum $X$, there is a metrizable compactum $Y$ and a co-elementary
map $f:X \to Y$.  In particular, every compactum is co-elementarily
equivalent to a metrizable one.\\  

Theorem \ref{5.3} has several sharper versions; one is Theorem 1.7 in
\cite{Ban14}, which sees the L\"{o}wenheim Theorem as a factorization
of maps.  The strongest version appears in \cite{Ban19}.

\subsection{Theorem.}\label{5.4}
(L\"{o}wenheim-Skolem Factorization Theorem, Theorem 3.1 in \cite{Ban19})
Let $f:X\to Y$ be a continuous surjection between compacta, with $\kappa$
an infinite cardinal such that $w(Y) \leq \kappa \leq w(X)$.  Then there
is a compactum $Z$ and continuous surjections $g: X \to Z$ and 
$h:Z \to Y$ such that $w(Z) = \kappa$, $g$ is a co-elementary map, and
$f = h\circ g$.\\  

\subsection{Remark.}\label{5.5}
When restricted to spaces in {\bf BS}, \ref{5.4} is an immediate
corollary of classical model theory.  In the absence of a Stone-like
duality, though, one must resort to other techniques.  The proof of
\ref{5.4} actually makes use of some Banach space theory.\\

Another line of inquiry regarding topological ultracoproducts concerns
the General Preservation Problem \ref{3.4}, with $\sum_{\cal D}X_i$ in
place of $\prod_{\cal D}X_i$.  In this new setting, we define a property
$P$ of compacta to be {\bf closed} if for any indexed family 
$\langle X_i: i \in I\rangle$ of compacta, and any ultrafilter $\cal D$
on $I$, $\sum_{\cal D}X_i$ has property $P$ whenever 
$\{i \in I: X_i\;\mbox{has property} \;P\} \in {\cal D}$.  $P$ is
{\bf open} if the complement of $P$ in {\bf CH} is closed. (Frequently we
speak of a subclass {\bf K} of {\bf CH} as being closed or open.)  

\subsection{Theorem.}\label{5.6}
The following properties of compacta are both closed and open:
\begin{quote}
``Being a continuum'' (Proposition 1.5 in \cite{Ban12}); ``being Boolean'' 
(Proposition
1.7 in \cite{Ban12}); ``having Lebesgue covering dimension $n$,
$n < \omega$,'' (essentially Theorem 2.2.2 in \cite{Ban12}); and
``being a decomposable/indecomposable continuum'' (Proposition 2.4.4 in
\cite{Ban12}, and Proposition 11 in \cite{Gur}).\\ 
\end{quote}

It follows from the above that 
the property of having infinite covering dimension is a closed property.
It is not an open property because $\sum_{\cal D}X_i$ will be infinite
dimensional as long as $\{i\in I: X_i\;\mbox{has dimension}\;n\}\in {\cal D}$
for arbitrarily large $n < \omega$.  The reader may be wondering whether
other dimension functions behave as well as covering dimension {\em vis \`{a}
vis\/} ultracoproducts, and the short answer is no:  There is a compactum
$X$, due to A. L. Lunc \cite{Nag, Pea} such that dim$(X) = 1$ and
ind$(X)$ = Ind$(X) = 2$ (where dim$(\;)$, ind$(\;)$ and Ind$(\;)$ are
covering dimension, small inductive dimension and large inductive dimension,
respectively).  Using \ref{5.3}, find a metrizable  
$Y \cong X$.  Then dim$(Y) = 1$ by
\ref{5.6}.  Since all three dimension functions agree for separable
metrizable spaces, we see that the two inductive dimension functions
are not preserved by co-elementary equivalence.

Recall that decomposability in a continuum $X$ means that $X$ is the 
union of two proper subcontinua; equivalently, it means that $X$ has
a proper subcontinuum with nonempty interior.  It is relatively easy
to show that the class of decomposable continua is closed; much less
trivial \cite{Gur} to show the same for the class of indecomposable
continua.  In \cite{Ban19} the class of {\bf $\kappa$-wide} compacta is
defined, for each cardinal $\kappa$.  Membership in this class amounts
to having a family of $\lambda$ pairwise disjoint proper subcontinua with 
nonempty
interiors, for each cardinal $\lambda < \kappa$; so decomposability
for a continuum is equivalent to being 1-wide, and all infinite locally
connected compacta are $\aleph_1$-wide.  Using a technique similar to
the one Gurevi\v{c} used to prove Proposition 11 in \cite{Gur}, one can
show that the class of $n$-wide compacta is both open and closed for
each $n<\omega$; consequently that any compactum co-elementarily equivalent
to a locally connected compactum is $\aleph_0$-wide.  The class of
$\aleph_0$-wide compacta is closed under co-elementary equivalence, but
this is hardly the case for the locally connected compacta. 

\subsection{Theorem.}\label{5.7}
(Corollary 14 in \cite{Gur}) Let $\cal D$ be a nonprincipal ultrafilter
on a countable set, with $X$ an infinite compactum.  Then 
$XI\backslash {\cal D}$ is not locally connected.\\

This result was used in \cite{Ban15} (along with regular ultrafilters and
the L\"{o}wenheim-Skolem Theorem) to obtain the following.

\subsection{Theorem.}\label{5.8}
(Theorem 2.10 in \cite{Ban15})  Let $\kappa$ be an infinite cardinal, and
$X$ an infinite compactum.  Then there is a compactum $Y$, of weight
$\kappa$, that is co-elementary equivalent to $X$, but not locally
connected.\\

The central role of local connectedness in the study of topological
ultracoproducts was discovered by R. Gurevi\v{c} in solving a problem I
raised in \cite{Ban12}.  In an exact analogy with the concept of
$\aleph_0$-categoricity in model theory, define a metrizable compactum $X$ 
to be {\bf $\aleph_0$-categorical} if there is no homeomorphically distinct
metrizable compactum that is co-elementarily equivalent to $X$.  For 
example, the Cantor discontinuum $2^{\omega}$ is $\aleph_0$-categorical
because its Boolean lattice of clopen sets is the unique (up to isomorphism)
countable atomless Boolean lattice, and the class of Boolean spaces is
both closed and open.  The problem I raised was whether the closed unit
interval $[0,1]$ (or {\em any} nontrivial metrizable continuum, for that 
matter) is $\aleph_0$-categorical, and Theorem \ref{5.8} provides a
negative answer.  (The same negative answer was given in \cite{Gur}, but
the proof of Proposition 15, a key step, was significantly incomplete.)
The question of the existence of $\aleph_0$-categorical continua remains
open, but we know from \ref{5.8} that any $\aleph_0$-categorical compactum
must fail to be locally connected.  (There is even more:  Using a Banach
version of the classic Ryll-Nardzewski Theorem from model theory, 
C. W. Henson \cite{Hen3} has informed me that $\aleph_0$-categorical
compacta must fail to be $\aleph_0$-wide.)

The concept of categoricity may be relativized to a subclass {\bf K}
of {\bf CH} in the obvious way.  Thus we could ask about the existence 
of metrizable compacta
in {\bf K} that are $\aleph_0$-categorical {\em relative to\/} {\bf K}.
When {\bf K} is the locally connected compacta, there is a satisfying
answer.  Recall that an {\bf arc} (resp {\bf simple closed curve}) is
a homeomorphic copy of the closed unit interval (resp. the standard unit
circle).

\subsection{Theorem.}\label{5.9}
(Theorem 0.6 in \cite{Ban13}) Arcs and simple closed curves
are $\aleph_0$-categorical relative to the class of locally connected
compacta.\\ 
 
Getting back to the General Preservation Problem \ref{3.4}, there is 
not much known about properties of a topological ultraproduct that
are conferred solely by the ultrafilter involved (in analogy with
\ref{3.8} and \ref{3.15}).  One such is due to K. Kunen \cite{Kun2}, and
uses a Banach space argument.

\subsection{Theorem.}\label{5.10}
(Kunen \cite{Kun2}) Let $\cal D$ be a regular ultrafilter on $I$, with
$X$ an infinite compactum.  Then $w(XI\backslash {\cal D}) =
w(X)^{|I|}$.\\ 

Recall that a $P_{\kappa}$-space is one for which intersections of fewer
than $\kappa$ open sets are open.  When $\kappa$ is uncountable, such
spaces are pseudofinite; hence infinite compacta can never be counted
among them.  There is a weakening of this property, however, that
compacta can subscribe to.  Call a space an {\bf almost-$P_{\kappa}$-space}
if nonempty intersections of fewer than $\kappa$ open sets have nonempty
interior.  The following is an easy consequence of the Additivity Lemma
\ref{3.8}, plus the fact that the topological ultracoproduct contains the
corresponding topological ultraproduct as a dense subspace.

\subsection{Theorem.}\label{5.11}
(Theorem 2.3.7 in \cite{Ban12}) If an ultrafilter is $\kappa$-regular, then
all topological ultracoproducts via that ultrafilter are 
almost-$P_{\kappa^+}$-spaces.\\

A little more significant is the following result about ultraco{\em powers\/}.
(Compare with \ref{3.15}.)

\subsection{Theorem.}\label{5.12}
(Theorem 2.3.17 in \cite{Ban12}) If an ultrafilter is $\kappa$-good, then all
topological ultracopowers via that ultrafilter are $B_{\kappa^+}$-spaces.\\

The rest of this section concerns what we have informally referred to as
the ``dualized model theory'' of compacta, in exact parallel (only with
the arrows reversed) with model-theoretic investigations of well-known
classes of relational structures (e.g., linear orders, graphs, groups,
fields, etc.).  As we saw above, the topological ultracoproduct allows for
the definition of co-elementary maps between compacta, as well as for the
creation of the hierarchy of classes of maps of level $\geq \alpha$ for 
any ordinal
$\alpha$.  When we restrict our attention to Boolean spaces, co-elementary
maps and maps of level $\geq \alpha$ are the Stone duals
of elementary embeddings and embeddings of level $\geq \alpha$, 
respectively, between Boolean lattices.  This basic correspondence
provides us with an abundance of facts about the Boolean setting that we
would like to extend to the compact Hausdorff setting.  Any failure of
extendability would give a new proof of the Banaschewski-Rosick\'{y}
Theorem \ref{5.1}; so far, however, there has been nothing but success
(or indecision).

The first obvious question that needs clearing up is whether the levels
really go beyond $\omega$, and the answer is no.

\subsection{Theorem.}\label{5.13}
(Hierarchy Theorem, Theorem 2.10 in \cite{Ban20}) Let $\alpha$ be any  
infinite ordinal.  Then the maps between compacta
that are of level $\geq \alpha$ are precisely the co-elementary maps.\\ 

This leads us to the second question, whether the composition of two maps
of level $\geq \alpha$ is also of level $\geq \alpha$.  As mentioned above,
the answer is yes, but a much stronger result is true.  The following is
a dualized version of (an easy generalization of) the Elementary Chains
Theorem of R. L. Vaught.  (I first conjectured the result in the mid
1970s, and it took twenty years to find a proof.)

\subsection{Theorem.}\label{5.14}
($\alpha$-Chains Theorem, Theorem 3.4 in \cite{Ban20}) 
Let $\langle X_n \stackrel{f_n}{\leftarrow} X_{n+1}: n < \omega \rangle$  
be a sequence of maps of level $\geq \alpha$ between compacta, with inverse
limit $X$ and limit maps $g_n:X \to X_n$, $n < \omega$.  Then each $g_n$
is a map of level $\geq \alpha$.\\
 
\subsection{Remark.}\label{5.15}
Without much ado, \ref{5.14} may be extended to arbitrary inverse systems
of compacta.\\

In any model-theoretic study of algebraic systems, the most commonly
investigated homomorphisms are the level $\geq 1$, or existential,
embeddings.  These are the ones arising from the classical study of
algebraically closed fields, for example.  When we look at the dual
notion of level $\geq 1$, or co-existential, maps between compacta, a
very rich theory emerges.  First of all, let us recall some
properties of compacta that are preserved by level $\geq 0$ maps
(alias continuous surjections).  These include: ``having 
cardinality (or weight) $\leq \kappa$ ($\kappa$ any cardinal);
``being connected;'' and ``being locally connected.''   
When we consider preservation by co-existential maps, we obtain
preservation for several important properties that are not generally 
preserved by continuous surjections.

\subsection{Theorem.}\label{5.16} 
(various results of \cite{Ban19})
The following properties are preserved
by co-existential maps:
\begin{quote}
``Being infinite;'' ``being disconnected;'' ``having covering dimension
$\leq n$'' ($n<\omega$);'' 
``being an indecomposable continuum;'' and ``being a hereditarily
indecomposable continuum.''\\ 
\end{quote}

\subsection{Remark.}\label{5.17}
Co-existential maps cannot raise covering dimension, but they can
lower it (Example 2.12 in \cite{Ban19}).  It is not hard to show 
that level $\geq 2$ maps between compacta must preserve covering
dimension.\\

An important tool in the proof of results such as \ref{5.16} is the 
following result, of interest in its own right.

\subsection{Theorem.}\label{5.18}
(Covering Lemma, Theorem 2.4 in \cite{Ban19}) Let $f:X \to Y$ be a
co-existential map between compacta.  Then there exists a
$\bigcup$-semilattice homomorphism $f^*$ from the subcompacta of $Y$
to the subcompacta of $X$ such that for each subcompactum $K$ of $Y$:
$(i)$ $f[f^*(K)] = K$; $(ii)$ $f^{-1}[U] \subseteq f^*(K)$ whenever $U$
is a $Y$-open set contained in $K$; $(iii)$ the restriction of $f$ to
$f^*(K)$ is a co-existential map from $f^*(K)$ to $K$; and $(iv)$ 
$f^*(K) \in \mbox{\bf K}$ whenever $K \in \mbox{\bf K}$ and $\mbox{\bf K}
\subseteq \mbox{\bf CH}$ is closed under ultracopowers and continuous
surjections.\\
 
An easy corollary of \ref{5.18} is the fact that co-existential maps between
compacta are {\bf weakly confluent}; i.e., possessed of the feature that
subcontinua of the range are themselves images of subcontinua of the 
domain.  If a subcontinuum of the range is the image of each component of
its pre-image, then the map is called {\bf confluent}.  Stronger still, a
continuous surjection is {\bf monotone} if pre-images of
subcontinua of the range are subcontinua of the domain.  

\subsection{Theorem.}\label{5.19}
(Theorem 2.7 in \cite{Ban19}) Let $f:X \to Y$ be a co-existential map
between compacta, where $Y$ is locally connected.  Then $f$ is monotone.\\

Theorem \ref{5.19} is a main ingredient for the following result; another
is Proposition 2.7 in \cite{Ban17}.

\subsection{Theorem.}\label{5.20}
Let $f:X\to Y$ be a function from an arc to a
compactum.  The following are equivalent:
\begin{quote}
$(i)$ $f$ is a co-existential map.\\
$(ii)$ $f$ is a co-elementary map.\\
$(iii)$ $Y$ is an arc and $f$ is a monotone continuous surjection.\\
\end{quote}

Recall that a class of relational structures is called {\bf elementary}
if it is the class of models of a first-order theory.  (This is the usage
in \cite{CK}.  In \cite{BS}, elementary classes are the classes of models
of a single sentence; what in \cite{CK} are called {\em basic\/} elementary
classes.)  From early work (1962) of T. E. Frayne, A. C. Morel and 
D. S. Scott (see Theorem 4.1.12 in \cite{CK}), a class is elementary if
and only if it is closed under the taking of ultraproducts and 
ultra{\em roots\/} (i.e., a structure is in the class if some ultrapower
of the structure is in the class).  This characterization is another
bridging theorem, allowing us to define a class $\mbox{\bf K} \subseteq
\mbox{\bf CH}$ to be {\bf co-elementary} if it is closed under the taking
of ultracoproducts and ultracoroots.  For example, all the classes 
(properties) mentioned in Theorem \ref{5.6} are co-elementary, since they
are both closed and open.  The class of compacta of infinite covering
dimension, while not being open, is still co-elementary.  The same may
be said for the class of $\aleph_0$-wide compacta (but certainly not
for the class of locally connected compacta, by \ref{5.7}).

An elementary class of relational structures is called {\bf model complete}
(see \cite{Maci})
if every embedding between members of that class is elementary.  Thus
we may define, in parallel fashion, the notion of {\bf model cocomplete
co-elementary class}.  (I apologize for so many {\em co}s.)  Because of
Stone Duality, plus the fact that the class of atomless Boolean lattices
is model complete, the class of self-dense Boolean spaces is a model
cocomplete class of compacta.
The following
is an exact analogue of Robinson's Test for model completeness, and
uses the $\omega$-Chains Theorem \ref{5.14}.

\subsection{Theorem.}\label{5.21}
(Robinson's Test, Theorem 5.1 in \cite{Ban19}) A co-elementary class of
compacta is model cocomplete if and only if every continuous surjection
between members of the class is a co-existential map.\\

In model theory, the Chang-\L o\'{s}-Suszko Theorem (see \cite{CK,Sho})
tells that an elementary class is the class of models of a set of 
universal-existential sentences if and only if the class is {\bf inductive};
i.e., closed under arbitrary chain unions.  In the compact Hausdorff
setting, we then define a co-elementary class to be {\bf co-inductive} if
that class is closed under inverse limits of chains of continuous 
surjections.  Examples of co-inductive co-elementary classes are 
{\bf CH}, {\bf BS} and {\bf CON} (the class of continua).
The co-elementary class of decomposable continua is not co-inductive; indeed
a favorite method of constructing indecomposable continua is to take
inverse limits of decomposable ones (see \cite{Nad}).  Define a
class {\bf K} of compacta to be {\bf $\kappa$-categorical}, where 
$\kappa$ is an infinite cardinal, if: $(i)$ {\bf K} contains compacta
of weight $\kappa$; and $(ii)$ any two members of {\bf K} of weight
$\kappa$ are homeomorphic.  The class of self-dense Boolean spaces, for
example, is $\aleph_0$-categorical.  The following is an exact analogue
of Lindstr\"{o}m's Test for model completeness, and uses Theorem \ref{5.21}
above, as well as a fair amount of topology.

\subsection{Theorem.}\label{5.22}
(Lindstr\"{o}m's Test, Theorem 6.4 in \cite{Ban19}) Any co-inductive
co-elementary class of compacta is model cocomplete, provided it contains
no finite members and is $\kappa$-categorical for some infinite cardinal
$\kappa$.\\

\subsection{Remark.}\label{5.23} Theorems \ref{5.21} and \ref{5.22} are
interesting and fairly hard to prove.  Unfortunately, they have proven
useless in finding interesting model cocomplete classes.  In particular,
we know of no model cocomplete classes of continua.\\

Model cocomplete co-elementary classes are interesting because, in some
sense, it is difficult to distinguish their members from one another.
This is especially true if they are also {\bf cocomplete}; i.e., consisting
of exactly one co-elementary equivalence class.  (It is not especially
hard to prove that every co-elementary equivalence class is closed, so
there is no problem finding cocomplete co-elementary classes.)  One way
to try to look for examples is via the study of {\em co-existential closure\/}.
Recall that in model theory, an $L$-structure $A$ is {\bf existentially
closed} relative to a class {\bf K} of $L$-structures, of which $A$ is
a member, if every embedding from $A$ into a member of {\bf K} is
existential.  Let $\mbox{\bf K}^e$ donote the members of {\bf K} that
are existentially closed relative to {\bf K}.  It is well known
(see \cite{CK}) that if {\bf K} is an inductive elementary class, then
each infinite $A \in \mbox{\bf K}$ embeds in some $A' \in \mbox{\bf K}^e$,
of cardinality $|A|$.  In certain special cases, $\mbox{\bf K}^e$ has a 
very elegant characterization.  For example, if {\bf K} is the class of
fields, then $\mbox{\bf K}^e$ is the class of algebraically closed
fields (Hilbert's {\em Nullstellensatz}).  Other examples include:
$(i)$ {\bf K} = the linear orderings without endpoints, 
$\mbox{\bf K}^e$ = the dense linear orderings without endpoints; 
$(ii)$ {\bf K} = the abelian groups, $\mbox{\bf K}^e$ = the divisible abelian
groups with infinitely many elements of each prime order.

We thus define a compactum $X \in \mbox{\bf K} \subseteq \mbox{\bf CH}$ to
be {\bf co-existentially closed} relative to {\bf K} if every continuous
surjection from a member of {\bf K} onto $X$ is co-existential.  Let 
$\mbox{\bf K}^c$ denote the members of {\bf K} that are co-existentially
closed relative to {\bf K}.  An exact analogue to the existence result
just cited is the following.

\subsection{Theorem.}\label{5.24}
(Level $\geq 1$ Existence, Theorem 6.1 in \cite{Ban19}) Let {\bf K} be
a co-inductive co-elementary class, with $X \in \mbox{\bf K}$ infinite.
Then $X$ is a continuous image of some $X' \in \mbox{\bf K}$, of 
weight $w(X)$.\\ 

Theorem \ref{5.24} applies, then, to the three co-inductive co-elementary
classes, {\bf CH}, {\bf BS} and {\bf CON} mentioned above.  The following
is not difficult to prove.

\subsection{Theorem.}\label{5.25}
(Proposition 6.2 in \cite{Ban19}) $\mbox{\bf CH}^c = \mbox{\bf BS}^c$ =
the self-dense Boolean spaces.\\

The nature of $\mbox{\bf CON}^c$ is apparently much more difficult to
discern.  If we can show it to be a co-elementary class, then, by
Robinson's Test \ref{5.21}, it is model-cocomplete.  (It is not hard to
show that $\mbox{\bf K}^c$ is closed under co-existential images when
{\bf K} is a co-elementary class.  Thus
to show $\mbox{\bf K}^c$ to be co-elementary, it suffices to show it is
closed under ultracoproducts.)  With a slight abuse of language, call
a member of $\mbox{\bf CON}^c$ a {\em co-existentially closed continuum\/}.
We know from \ref{5.24} that co-existentially closed continua abound,
but the process used to construct them involves direct limits of lattices,
and is not very informative.  We have very few criteria to decide whether
a given continuum is co-existentially closed; what we know so far is
the following.

\subsection{Theorem.}\label{5.26}
(Theorem 4.5 in \cite{Ban20}) Every co-existentially closed continuum
is an indecomposable continuum of covering dimension one.\\

\section{Related Constructions.}\label{6}
Starting with an $I$-indexed family $\langle X_i: i \in I\rangle$
of topological spaces, the box product topology on the set
$\prod_{i \in I}X_i$ is defined by declaring the {\em open boxes\/}
$\prod_{i \in I}U_i$ as basic open sets, where the sets $U_i$ are open
subsets of $X_i$, $i \in I$.  Alternatively, one forms the Tychonov
product topology by restricting attention to those open boxes having
the property that $\{i\in I: U_i \neq X_i\}$ is finite.  In \cite{Kni},
C. J. Knight combines these two formations under a common generalization,
the {\bf $\cal I$-product topologies} for any ideal ${\cal I}$ of subsets of $I$
(so $\emptyset \in {\cal I}$, and $\cal I$ is closed under subsets and
finite unions) as follows:  Take as open base all open boxes
$\prod_{i \in I}U_i$ such that $\{i \in I: U_i \neq X_i\} \in {\cal I}$.
Then the box (resp., Tychonov) product topology  
is the $\cal I$-product topology for  
${\cal I} := {\wp}(I)$ (resp., ${\cal I}:= 
\{J \subseteq I: J\;\mbox{finite}\;\}$).  (For the trivial ideal
${\cal I} := \{\emptyset\}$, one trivially obtains the trivial topology.)
The collective name for these $\cal I$-product formations, for various
ideals $\cal I$, is known as the {\bf ideal product topology}.

In \cite{GK}, M. Z. Grulovi\'{c} and M. S. Kurili\'{c} add a new
ingredient to the pot, creating a further generalization that now
takes in
all ideal product topologies, as well as all reduced product topologies.
Known as the {\bf reduced ideal product topology}, it comprises the
{\bf ${\cal FI}$-product topologies} for any pair 
$\langle {\cal F},{\cal I}\rangle$,
where $\cal F$ (resp., $\cal I$) is a filter (resp., an ideal) on $I$:
First one takes the $\cal I$-product topology on $\prod_{i \in I}X_i$;
then forms the obvious quotient topology on the reduced product
$\prod_{\cal F}X_i$ of underlying sets.  Denote this new space by
$\prod_{\cal F}^{\cal I}X_i$.  Then the topological reduced product
$\prod_{\cal F}X_i$ of \S 3 is   
$\prod_{\cal F}^{{\wp}(I)}X_i$ in this notation.  Also, when $\cal F$ includes
all the complements of members of $\cal I$, it follows that 
$\prod_{\cal F}^{\cal I}X_i$ has the trivial topology.

Define a filter-ideal pair  $\langle {\cal F},{\cal I}\rangle$ on $I$
to satisfy the {\bf density condition} if for every $A \in {\cal F}$
and every $B \notin {\cal F}$, there exists a $C \in {\cal I}$ such
that $C \subseteq A\backslash B$ and $I\backslash C \notin {\cal F}$.
(The use of the word {\em density\/} in this definition is justified by
the following observation.
Consider the quotient partially ordered
set ${\wp}(I)/{\cal F}$, where $A,B \subseteq I$ are identified if
$A \cap F = B \cap F$ for some $F \in {\cal F}$.
Then the density condition amounts to the condition
that every nonbottom element of
${\wp}(I)/{\cal F}$ dominates a nonbottom element of 
${\cal I}/{\cal F}$.) Note that $\langle {\cal F},{\wp}(I)\rangle$ 
satisfies the density condition when $\cal F$ is a proper
filter (given $A$ and $B$, just let $C$ be $A\backslash B$), and that
$\langle \{I\}, {\cal I}\rangle$ satisfies the density condition when
every nonempty subset of $I$ contains a nonempty member of $\cal I$.  Also 
note that if $\cal F$ includes all the complements of members of ${\cal I}$, 
then $\langle {\cal F},{\cal I}\rangle$ does not
satisfy the density condition.

The main contribution of \cite{GK} is to connect the density condition
with the preservation of the separation axioms by reduced ideal products
(in a manner not entirely unlike the style of Theorem \ref{4.1}).
For a topological property $P$, say that a filter-ideal pair
$\langle {\cal F},{\cal I}\rangle$ {\bf preserves} $P$ if for any
$I$-indexed family $\langle X_i:i \in I\rangle$, 
$\prod_{\cal F}^{\cal I}X_i$ has property $P$ whenever
$\{i \in I: X_i\;\mbox{has property}\;P\} \in {\cal F}$. 

\subsection{Theorem.}\label{6.1}
(Grulovi\'{c}-Kurili\'{c} \cite{GK}) Let property $P$ be any of the
separation axioms $T_r$, $r \in \{0,1,2,3,3.5\}$.  Then a filter-ideal
pair $\langle {\cal F},{\cal I}\rangle$ preserves $P$ if and only if 
it satisfies the density condition.\\

\section{Open Problems.}\label{7}
Here are some of the problems that have remained tantalizingly open.

\subsection{}\label{7.1} 
(See \ref{3.6}.) Can a topological ultraproduct be normal without being
paracompact?\\

\subsection{}\label{7.2} 
(See \ref{3.15}, \ref{3.16}.) If all topological ultraproducts via $\cal D$ are
$\kappa$-Baire, as well as $\lambda^+$-additive for all $\lambda < \kappa$,
is $\cal D$ necessarily $\kappa$-good?\\

\subsection{}\label{7.3} 
(See \ref{3.17}.) Is there a nice topological characterization of 
${\Bbb Q}_{\alpha}$ for uncountable $\alpha$? (Candidate: being
regular self-dense, of cardinality = weight $= \alpha$, 
$\alpha$-additive and $\alpha$-Baire.  It is definitely not enough to
exclude the ``$\alpha$-Baire'' part, as Example 3.11 in \cite{Ban9}
shows.)\\

\subsection{}\label{7.4} 
(See \ref{3.18}, \ref{3.19}.) Is pseudofiniteness (i.e., anticompactness) 
preserved by topological ultraproducts?\\

\subsection{}\label{7.5}
(See \ref{4.1}, \ref{4.2}.) Are there (consistently) any $T_4$-ultrafilters?\\

\subsection{}\label{7.6}
(A representation problem.) 
If $X$ and $Y$ are co-elementarily equivalent compacta, can one always find
lattice bases ${\cal A}$ for $X$ and ${\cal B}$ for $Y$ such that 
${\cal A} \equiv {\cal B}$?\\

\subsection{}\label{7.7}
(See \ref{5.8} and subsequent discussion.) Do nontrivial
$\aleph_0$-categorical continua exist?\\

\subsection{}\label{7.8}
(See \ref{5.9}.)  What Peano continua are there (besides arcs
and simple closed curves) that are $\aleph_0$-categorical relative to
the class of locally connected compacta?\\ 

\subsection{}\label{7.9}
(See \ref{3.15}, \ref{5.12}.) Do $\kappa$-good ultrafilters
create ultracoproduct compacta that are $B_{\kappa^+}$-spaces?\\

\subsection{}\label{7.10}
(See \ref{5.18}, \ref{5.19}.) Are co-existential maps always
confluent?\\
 
\subsection{}\label{7.11}
(See \ref{5.18}.) Is there a true version of the Covering Lemma
where {\em co-existential\/} is replaced by {\em of level $\geq \alpha$\/}
($\alpha \geq 1$)?\\
 
\subsection{}\label{7.12}
(See \ref{5.23}.) Are there any model cocomplete co-elementary classes
that contain nontrivial continua?\\
 
\subsection{}\label{7.13}
(See \ref{5.26}.) Is the class of co-existentially closed continua a
co-elementary class?\\

\subsection{}\label{7.14}
(See \ref{5.26}.) Are any of the familiar one-dimensional indecomposable
continua (e.g., the solenoids, the bucket handle continua, the lakes
of Wada) co-existentially closed?\\

\end{document}